\documentclass[a4paper,12pt]{article}

\usepackage[utf8]{inputenc}
\usepackage[T1]{fontenc}
\usepackage[left=2.5cm,right=2.5cm,top=1.5cm,bottom=2cm]{geometry}
\usepackage{amsmath,amssymb,makeidx}
\usepackage[english]{babel}

\usepackage[language=english,backend=biber]{biblatex}
\usepackage{fourier}
\usepackage{graphics}
\usepackage{graphicx}
\usepackage{amsthm}
\usepackage{faktor}
\usepackage{thmtools}
\usepackage[hyperindex=true]{hyperref}
\usepackage{color}
\usepackage{enumitem}
\newlist{itemizeth}{itemize}{2}
\setlist[itemizeth]{label=\textbullet,noitemsep,topsep=0 mm}
\newlist{enumerateth}{enumerate}{2}
\setlist[enumerateth]{label*=\textbf{\alph*)},itemsep=0.2mm}
\setlist[itemize]{fullwidth,topsep=0 mm, leftmargin=0pt}
\setlist[enumerate]{label*=\textbf{\arabic*)},itemsep=0.2mm, leftmargin=0pt}
\declaretheoremstyle[
title=Proof,
numbered=no,
headfont=\normalfont\bfseries,
notefont=\bfseries, notebraces={}{},
postheadspace=17pt,
bodyfont=\normalfont,
headindent=0pt,
qed=\qedsymbol,
spacebelow=3,
]{demostyle}

\declaretheoremstyle[
headfont=\normalfont\bfseries,
notefont=\mdseries, notebraces={(}{)},
bodyfont=\normalfont,
thmbox=S
]{standard}

\declaretheoremstyle[
title=Theorem,
headfont=\normalfont\scshape,
notefont=\mdseries, notebraces={(}{)},
thmbox=S,
]{theo}

\declaretheoremstyle[
title=Theorem,
headfont=\scshape\bfseries,
notefont=\mdseries, notebraces={(}{)},
thmbox=L,
]{import}

\declaretheoremstyle[headfont=\normalfont\itshape, 
numbered=no
]{rem}

\declaretheorem[title=Definition,style=standard]{defi}
\declaretheorem[title=Proposition,style=standard]{prop}
\declaretheorem[title=Lemma,style=standard]{lem}

\declaretheorem[style=theo]{Th}

\declaretheorem[title=Remark,style=rem]{rqu}
\declaretheorem[title=Remarks,style=rem]{rqus}
\declaretheorem[style=demostyle]{dem}

\makeatletter
\newcommand*{\house}[1]{%
   \mathord{%
     \mathpalette\@house{#1}%
   }%
}
\newcommand*{\@house}[2]{%
   \dimen@=\fontdimen8 %
       \ifx#1\scriptscriptstyle\scriptscriptfont
       \else\ifx#1\scriptstyle\scriptfont
       \else\textfont\fi\fi
       3 %
   \sbox0{%
     $#1%
       \vrule width\dimen@\relax
       \overline{%
         \kern2\dimen@
         \begingroup 
           #2%
         \endgroup
         \kern2\dimen@
       }%
       \vrule width\dimen@\relax
       \mathsurround=1.5\dimen@ 
     $%
   }%
   \ht0=\dimexpr\ht0-\dimen@\relax
   \dp0=\dimexpr\dp0+2\dimen@\relax
   \vbox{%
     \kern\dimen@ 
     \copy0 %
   }%
}

\newcommand{\R}{\mathbb{R}}
\newcommand{\N}{\mathbb{N}}

\newcommand{\Z}{\mathbb{Z}}
\newcommand{\Q}{\mathbb{Q}}
\newcommand{\C}{\mathbb{C}}
\newcommand{\K}{\mathbb{K}}
\newcommand{\Oal}{\mathcal{O}}
\newcommand{\Qbar}{\overline{\mathbb{Q}}}
\newcommand{\Spec}{\mathrm{Spec} \,}

\newcommand{\GL}{\mathrm{GL}}
\newcommand{\Pro}{\mathbb{P}}
\newcommand{\ord}{\mathrm{ord} \,}
\newcommand{\ssi}{\Leftrightarrow}
\newcommand{\gf}{\dfrac}
\newcommand{\Vect}{\mathrm{Span}}

\newcommand{\den}{\mathrm{den}}
\newcommand{\ddz}{\gf{\mathrm{d}}{\mathrm{d}z}}
\newcommand{\singze}{\left\lbrace 0 \right\rbrace}
\newcommand{\db}{\mathcal{D}}
\newcommand{\fonction}[5]{\begin{array}[t]{lrcl} 
#1: & #2 & \longrightarrow & #3 \\
    & #4 & \longmapsto & #5 \end{array}}

\title{On the linear independence of values of $G$-functions}
\date{\today}
\author{Gabriel Lepetit}
\makeindex
\begin{document}

\maketitle

\begin{abstract}
We consider a $G$-function $F(z)=\sum_{k=0}^{\infty} A_k z^k \in \K\llbracket z \rrbracket$, where $\K$ is a number field, of radius of convergence $R$ and annihilated by the $G$-operator $L \in \K(z)[\mathrm{d}/\mathrm{d}z]$, and a parameter $\beta \in \Q \setminus \Z_{\leqslant 0}$. We define a family of $G$-functions $F_{\beta,n}^{[s]}(z)=\sum_{k=0}^{\infty} \gf{A_k}{(k+\beta+n)^s} z^{k+n}$ indexed by the integers $s$ and $n$. Fix $\alpha \in \K^* \cap D(0,R)$. Let $\Phi_{\alpha,\beta,S}$ be the $\K$-vector space generated by the values $F_{\beta,n}^{[s]}(\alpha)$, $n \in \N$, $0 \leqslant s \leqslant S$. We show that there exist some positive constants $u_{\K,F,\beta}$ and $v_{F,\beta}$ such that $u_{\K,F,\beta} \log(S) \leqslant \dim_{\K} \Phi_{\alpha,\beta,S} \leqslant v_{F,\beta} S$.
This generalizes a previous theorem of Fischler and Rivoal (2017), which is the case $\beta=0$. Our proof is an adaptation of their article \cite{FRivoal}, making use of the André-Chudnovsky-Katz Theorem on the structure of the $G$-operators and of the saddle point method.
Moreover, we rely on Dwork and André's quantitative results on the size of $G$-operators to obtain an explicit formula for the constant $u_{\K,F,\beta}$, which was not given in \cite{FRivoal} in the case $\beta=0$.
\end{abstract}
\bigskip

\section{Introduction}\label{sec:intro}

A \emph{$G$-function} is a power series $f(z)=\sum\limits_{k=0}^{\infty} a_k z^k \in \Qbar\llbracket z \rrbracket$ satisfying the three following assumptions:

\begin{enumerateth}
\item $f$ is solution of a nonzero linear differential equation with coefficients in $\Qbar(z)$;
\item There exists  $C_1 >0$ such that $\forall k \in \N, \; \house{a_k} \leqslant C_1^{k+1}$, where $\house{a_k}$ is the \emph{house} of $a_k$, \emph{i.e.} the maximum of the absolute values of the Galois conjugates of $a_k$;
\item There exists $C_2 >0$ such that $\forall k \in \N, \;  \mathrm{den}(a_0, \dots,a_k) \leqslant C_2^{k+1}$, where $\mathrm{den}(a_0, \dots,a_k)$ is the \emph{denominator} of $a_0, \dots, a_k$, \emph{i.e.} the smallest integer $d \in \N^*$ such that $d a_0, \dots, d a_k$ are algebraic integers.
\end{enumerateth}

This family of special functions has been studied together with the family of $E$-functions, which are the functions $f(z)=\sum\limits_{k=0}^{\infty} \gf{a_k}{k!} z^k$ satisfying \textbf{a)} and such that the $a_k$ satisfy the conditions \textbf{b)} and \textbf{c)}, since Siegel defined them in 1929 \cite{Siegelarticle}. The most basic example of $G$-function, which gives it its name, is the geometric series $f(z)=-\sum\limits_{k=0}^{\infty} z^k=1/(1-z)$. Other examples include $$\log(1-z)=\sum_{k=1}^{\infty} \gf{z^k}{k}, \quad \mathrm{Li}_s(z)=\sum\limits_{k=1}^{\infty} \gf{z^k}{k^s}, \quad \arctan(z)=\sum_{k=0}^{\infty} (-1)^k \gf{z^{2k+1}}{2k+1}$$ and the family of hypergeometric functions with rational parameters: if $\boldsymbol{\alpha}=(\alpha_1, \dots, \alpha_{\nu}) \in \Q^{\nu}$ and $\boldsymbol{\beta}=(\beta_1, \dots, \beta_{\nu-1}) \in (\Q \setminus \Z_{\leqslant 0})^{\nu-1}$, $$_{\nu} F_{\nu-1} (\boldsymbol{\alpha} ; \boldsymbol{\beta} ;z):=\sum\limits_{k=0}^{\infty} \gf{(\alpha_1)_k \dots (\alpha_{\nu})_k}{(\beta_1)_k \dots (\beta_{\nu-1})_k k!} z^k, \quad $$  where for $x \in \C$ and $k \geqslant 1$, $(x)_k :=x(x+1)(x+2) \dots (x+k-1)$, $(x)_0=1$, is the \emph{Pochhammer symbol}.

In this paper, we are going to rely on the theory of $G$-functions developed by André, Bombieri, Chudnovsky, Katz and others. Its main result can be synthetised as follows: the minimal nonzero linear differential equation on $\Qbar(z)$ associated with a $G$-function belong to a specific class of differential operators, called \emph{$G$-operators}. Every $G$-operator of order $\mu$ is Fuchsian and admits a basis of solutions around every point $a$ of $\Pro^{1}(\Qbar)$ of the form $(f_1(z-a), \dots, f_{\mu}(z-a)) (z-a)^{C_u}$, where the $f_i(u)$ are $G$-functions and $C_u \in \mathcal{M}_{\mu}(\Q)$. See \cite{Andre} or \cite{Dwork} for an extensive review of the theory of $G$-functions.

\bigskip

Our goal is to study the following problem, first considered in a special case (\emph{i.e.}, $\beta=0$) by Fischler and Rivoal (\cite{FRivoal}), involving $G$-functions and $G$-operators.
Let $\K$ a number field and $F(z)=\sum\limits_{k=0}^{\infty} A_k z^k \in \K\llbracket z \rrbracket$ a non polynomial $G$-function. Let $L \in \Qbar\left[z, \mathrm{d}/\mathrm{d}z \right] \setminus\singze$ be an operator such that $L(F(z))=0$ and of minimal order $\mu$ for $F$.

Take a parameter $\beta \in \Q \setminus \Z_{\leqslant 0}$, that will remain fixed in the rest of the paper. For $n \in \N^*$ and $s \in \N$, we define the $G$-functions $$F_{\beta,n}^{[s]}(z)=\sum\limits_{k=0}^{\infty} \gf{A_k}{(k+\beta+n)^s} z^{k+n}.$$ These are related to iterated primitives of $F(z)$. The aim of this article is to find upper and lower bounds of the dimension of $$\Phi_{\alpha,\beta, S}:=\Vect_{\K}\left(F_{\beta,n}^{[s]}(\alpha), \; n \in \N^*, \;\; 0 \leqslant s \leqslant S\right)$$ when $S$ is a large enough integer and $\alpha \in \K$, $0<|\alpha| <R$ with $R$ the radius of convergence of~$F$. Note that it is not obvious that $\Phi_{\alpha,\beta,S}$ has finite dimension. Precisely, we want to prove the following theorem.

\begin{Th} \label{th:chapitre2}
Assume that $F$ is not a polynomial. Then for $S$ large enough, the following inequality holds: $$ \gf{1+o(1)}{[\K:\Q] C(F,\beta)} \log(S) \leqslant \dim_{\K} \Phi_{\alpha,\beta,S} \leqslant \ell_0(\beta) S+\mu.$$ Here, if $\delta=\deg_z(L)$ and $\omega$ is the order of $0$ as a singularity of $L$, $\ell_0(\beta)$ is defined as the maximum of $\ell:=\delta-\omega$ and the numbers $f-\beta$ when $f$ runs through the exponents of $L$ at infinity such that $f-\beta \in \N$, and $C(F,\beta)$ is a positive constant depending only on $F$ and $\beta$, and not $\alpha$.
\end{Th}

If $F(z) \in \K[z]$, then $F_{\beta,n}^{[s]}(z) \in \K[z]$ as well and $\Phi_{\alpha,\beta,S} \subset \K$.

In \cite{FRivoal}, Fischler and Rivoal proved this theorem with $\beta=0$. Their goal was to generalize previous results of Rivoal (\cite{Rivoal2003}, for $\alpha \in \Q$) and Marcovecchio \cite{Marcovecchio} on the dimension of the $\K$-vector space spanned by  $\mathrm{Li}_s(\alpha)$, $0 \leqslant s \leqslant S$, for $\alpha \in \K$, $0 < |\alpha| <1$, where $\mathrm{Li}_s(z):=\sum\limits_{k=1}^{\infty} \gf{z^k}{k^s}$ is the $s$-th polylogarithm function. Indeed, if we set $A_k=1$ for every positive integer $k$, the family of functions $\left(F_{0,n}^{[s]}(z)\right)_{n,s}$ is the family of the polylogarithms $\mathrm{Li}_s(z)$ up to an additive polynomial term. Using a different method based on a generalization of Shidlovskii's lemma,  Fischler and Rivoal later proved in \cite{FRivoalII} that Theorem \ref{th:chapitre2} was also true for $\beta=0$ and $\alpha$ in a domain that is star-shaped  at 0 in which the open disc of convergence of $F$ is strictly contained. We don't know if this also holds for any rational $\beta$, and it seems to be a difficult task.

We are going to adapt their approach in \cite{FRivoal} to the more general case we are interested in. In a first part, we rely on the properties of $G$-function of $F$ to find a recurrence relation between the functions $F_{\beta,n}^{[s]}(z)$, which will prove the upper bound of the theorem. In a second part, we will study the asymptotic behavior as $n \rightarrow +\infty$ of a power series $T_{S,r,n}(z)$, which is a linear form in the $F_{\beta,n}^{[s]}(z)$, in order to use a linear independence criterion à la Nesterenko due to Fischler and Rivoal, leading to the lower bound of Theorem \ref{th:chapitre2}. The key tool for this will be the saddle point method. 

In Section \ref{sec:explicitC}, we will give an original explicit expression of the constant $C(F,\beta)$. To this end, we recall in Section \ref{sec:sizeGop} results of Dwork \cite{Dwork}, André \cite{Andre} and of the author \cite{LepetitSize} on the notion of \emph{size} of a $G$-operator, encoding a condition of moderate growth on some denominators, the \emph{Galochkin condition}. In particular, an explicit version of Chudnovsky's Theorem gives a relation between the size of a $G$-function, encoding the conditions \textbf{b)} and \textbf{c)} of the definition above, and the size of its minimal operator.

After simplification, the estimation of $C(F,\beta)$ we obtain ultimately depends on $\beta$ (in fact, on its denominator), on arithmetic and analytic invariants of the minimal operator $L$ of $F$ and on the size of $F$ itself. 

In order to compute $C(F,\beta)$, it is possible and more convenient to rewrite $L$ in the form \begin{equation}\label{eq:L=sommeQj} L= \gf{z^{\omega-\mu}}{u} \sum_{j=0}^{\ell} z^j Q_j(\theta+j)\;, \quad \theta=z \gf{\mathrm{d}}{\mathrm{d}z}\end{equation} with $\omega \in \N^*$, $Q_j(X) \in \Oal_{\K}[X]$ and $u \in \N^*$. We notice that if $\ell=0$, then the only power series solutions of $L(y(z))=0$ are polynomials (see the remark after Equation \eqref{eq:recurrenceLbeta}, Section \ref{sec:explicitC}).

We will prove the following theorem:

\begin{Th} \label{th:calculC1etC2}
Assume that $F$ is not a polynomial. We denote by $\db$ the denominator of $\beta$. Then the integer $\ell$ defined in \eqref{eq:L=sommeQj} below is $\geqslant 1$ and a suitable constant $C(F,\beta)$ in Theorem \ref{th:chapitre2} is \begin{equation} \label{eq:expressionCFbeta}
    C(F,\beta)=\log(2eC_1(F) C_2(F,\beta))
\end{equation} with \begin{equation}\label{eq:expressionC1} C_1(F):=\max\left(1,\house{\gamma_0/\gamma_{\ell}}, \Phi_0(L)^{\max(1,\ell-1)}\right)\end{equation} and \begin{multline}\label{eq:expressionC2} C_2(F,\beta) := \\ \den\left(1/\gamma_0\right)^3 \den(\mathbf{e}, \beta)^{6 \mu} \exp\left(3 \max(1,\ell-1)[\K:\Q] \Lambda_0(L,\beta)+3(\mu+1)\den(\mathbf{f}, \beta) \right),\end{multline}
where the polynomials $Q_0(X)=\gamma_0 \prod\limits_{i=1}^{\mu} (X-e_i)$ and $Q_{\ell}(X)=\gamma_{\ell} \prod\limits_{i=1}^{\mu} (X+f_i-\ell)$ are defined in~\eqref{eq:L=sommeQj}. 
\end{Th}

The numbers $e_i$ (resp. $f_j$) are congruent modulo $\Z$ to the exponents of $L$ at $0$ (resp. $\infty$) and are therefore rational numbers by Katz's Theorem (\cite[p. 98]{Dwork}), since $L$ annihilates the $G$-function $F$. The numbers $\Lambda_0(L,\beta)$ and $\Phi_0(L)$ will be defined respectively in formulas \eqref{eq:defLambda0} and \eqref{eq:definitionPhi0} of Section \ref{sec:explicitC}. This theorem provides a constant $C(F,\beta)$ that eventually only depends on $F$ and the denominator $\db$ of $\beta$.

The terms $C_1(F)$ and $C_2(F,\beta)$ making up the constant $C(F,\beta)$ arise from very different computations: $C_1$ can be seen as the "analytic" part of $C$ whereas $C_2$ is related to arithmetic invariants of $F$, $L$ and $\beta$.

\medskip

In Section \ref{sec:examples}, we will end the paper by making explicit the results of Theorems \ref{th:chapitre2} and \ref{th:calculC1etC2} in the cases of classical examples, including polylogarithms, hypergeometric functions, and the generating function of the Apéry numbers.

\bigskip

\textbf{Acknowledgements:} I thank T. Rivoal for carefully reading this paper and for his useful comments and remarks that improved it substantially. I also thank the anonymous referee for pointing out some imprecisions and mistakes in the manuscript.

\section{Proof of the main result}\label{sec:cas1D}

\subsection{A recurrence relation between the $F_{\beta,n}^{[s]}(z)$}\label{subsec:recurrenceFns}

As $F$ is a $G$-function, the nonzero minimal operator $L$ of $F$ is a $G$-operator by Chudnovsky's Theorem (see \cite[p. 267]{Dwork}). In particular, the André-Chudnovsky-Katz Theorem (cf \cite[p.719]{AndregevreyI}) mentioned in Section \ref{sec:intro} states that $L$ is a Fuchsian operator with rational exponents at every point of $\Pro^1(\Qbar)$.
Relying on this property, we are going to obtain a recurrence relation between the series $F_{\beta,n}^{[s]}(z)$, when $s \in \N$ and  $n \in \N^*$ (Proposition \ref{prop:recurrenceFns} below). Here, and in all that follows, $\N$ (resp. $\N^*$) denotes the set of non-negative (resp. positive) integers. 

This algebraic method will be the key argument to obtain the upper bound in Theorem \ref{th:chapitre2} and will also be useful in Subsection \ref{subsec:auxiliaryseries}, where we will prove the lower bound.

By \cite[Lemma 1, p.  11]{FRivoal}, there exist some polynomials $Q_j(X) \in \Oal_{\K}[X]$ and $u \in \N^*$ such that  \begin{equation}\label{eq:reecritureL} u z^{\mu-\omega} L=\sum_{j=0}^{\ell} z^j Q_j(\theta+j),\end{equation} with $\theta=z \mathrm{d}/\mathrm{d}z$, $\mu$ the order of $L$, $\omega$ the multiplicity of $0$ as a singularity of $L$ and $\ell=\delta-\omega$ where $\delta$ is the degree in $z$ of $L$.

\begin{lem} \label{lem:Lbetaopmin}
Define, for $j \in \{ 0, \dots, \ell \}$, $Q_{j, \beta}(X)=\db^{\mu} Q_j(X-\beta)$, where $\db=\den(\beta)$. The operator $$L_{\beta}=\sum\limits_{j=0}^{\ell} z^j Q_{j,\beta}(\theta+j)$$ is an operator in $\Qbar\left[z,\mathrm{d}/\mathrm{d}z\right] \setminus \{0\}$ of minimal order for $z^{\beta} F(z)$.
\end{lem}

Before proving Lemma 1, we mention the following consequence of Chudnovsky's Theorem stated by Dwork (\cite[Corollary 4.2, p. 299]{Dwork}). 

\begin{prop}\label{prop:Gopbasesol}
Let $L$ be an operator in $\Qbar(z)\left[\mathrm{d}/\mathrm{d}z\right]$ such that the differential equation $L(y(z))=0$ has a basis of solutions around $0$ of the form $(f_1(z), \dots, f_{\mu}(z)) z^A$, where the $f_i(z)$ are $G$-functions and the matrix $A \in \mathcal{M}_{\mu}(\Qbar)$ has rational eigenvalues. Then $L$ is a $G$-operator.
\end{prop}

We recall that $z^A$ is defined in a simply connected open subset $\Omega$ of $\C^*$, as $z^A :=\exp(A \log(z))=\sum\limits_{k=0}^{\infty} \log(z)^k A^k / k! $ for $z \in \Omega$, where $\log$ is a determination of the complex logarithm on $\Omega$.

Proposition \ref{prop:Gopbasesol} implies that $L_{\beta}$ is a $G$-operator. Indeed, if we set a basis of solutions of $L(y(z))=0$ around $0$ of the form $(f_1(z), \dots, f_{\mu}(z)) z^{C}$, where $C \in \mathcal{M}_{\mu}(\Qbar)$ has eigenvalues in $\Q$, a basis of solutions of $L_{\beta}(y(z))=0$ around $0$ is $(f_1(z), \dots, f_{\mu}(z)) z^{C+\beta I_{\mu}}$.

\begin{dem}[of Lemma \ref{lem:Lbetaopmin}]
We begin by the following observation: for all $m, j \in \N$, $$(\theta-\beta+j)^m (z^{\beta} F(z))=z^{\beta} (\theta+j)^m(F(z)).$$  It is enough to prove it for $F(z)=z^k$. In that case, for $m=1$, $$(\theta-\beta+j)(z^\beta z^k)=(\beta+k)z^{\beta+k}+(j-\beta)z^{\beta+k}=z^{\beta}(k+j)z^k=z^{\beta}(\theta+j) z^k$$ and the result follows by induction on m. 

Now, with $Q_j=\sum\limits_{m=0}^{d_j} \rho_{j,m} X^m$, we have \begin{align*}
L_{\beta}(z^{\beta} F(z)) &= \db^{\mu} \sum_{j=0}^{\ell} \sum_{m=0}^{d_j} z^j \rho_{j,m} (\theta-\beta+j)^m (z^{\beta} F(z))=\db^{\mu} z^{\beta} \sum_{j=0}^{\ell} \sum_{m=0}^{d_j} z^j \rho_{j,m} (\theta+j)^m (F(z)) \\
&= \db^{\mu} z^{\beta} \sum_{j=0}^{\ell} z^j Q_j(\theta+j)(F(z))=0.
\end{align*} 

We note that $L_{\beta}$ has the same order as $L$. Let us now prove that this order is the minimal one for $z^{\beta} F(z)$. Let $\widetilde{L}$ be an operator of minimal order for $\widetilde{F}(z):=z^{\beta} F(z)$. Then $F(z)=z^{-\beta} \widetilde{F}(z)$, so $\widetilde{L}_{-\beta}(F)=0$. Thus $\ord(L) \leqslant \ord(\widetilde{L}_{-\beta})=\ord(\widetilde{L})$, since $L$ is minimal for $F$. On the other hand, $L_{\beta}(\widetilde{F})=0$, so that the minimality of $\widetilde{L}$ yields $\ord(\widetilde{L}) \leqslant \ord(L)$. Finally, $\widetilde{L}$ has the same order as $L$ and $L_{\beta}$, so $L_{\beta}$ is indeed minimal.
\end{dem}

In a similar way as in \cite[Lemma 3, p. 17]{FRivoal}, we obtain the following key lemma:

\begin{lem} \label{lemmerecurrenceFns}
For any fixed $s \in \N^*$, the sequence of functions $(F_{\beta,n}^{[s]}(z))_{n \geqslant 1}$ satisfies the following inhomogeneous recurrence relation:
$$\forall n \geqslant 1, \quad \sum_{j=0}^\ell Q_{j, \beta}(-n) F_{\beta, n+j}^{[s]}(z)=\sum_{j=0}^\ell \sum\limits_{t=1}^{s-1} \gamma_{j,n,t,s,\beta} F_{\beta, n+j}^{[t]}(z)+\sum\limits_{j=0}^\ell z^{n+j} B_{j,n,s,\beta}(\theta) F(z) $$ where $\gamma_{j,n,t,s,\beta} \in \Oal_\K$ and each polynomial $B_{j,n,s,\beta}(X) \in \Oal_\K[X]$ has degree at most $d_j-s$.\end{lem}

\begin{dem}
Let us proceed by induction on $s \geqslant 1$. 

\begin{itemize}
\item For $s=1$, let us remark that for $u \in \N$, $$\int_{0}^z x^{\beta+u} F(x) \mathrm{d}x=\sum_{k=0}^{\infty} A_k \int_{0}^z x^{\beta+u+k} \mathrm{d}x = \sum_{k=0}^{\infty} \gf{A_k}{\beta+u+k+1} z^{\beta+u+k+1}=z^{\beta} F_{\beta,u+1}^{[1]}(z).$$ Hence, if we set $L_1=u z^{\mu-\omega} L$ as in \eqref{eq:reecritureL} above, we have $$0 = \int_{0}^z x^{\beta+n-1} L_1(F(x)) \mathrm{d}x=\sum_{j=0}^\ell \sum_{m=0}^{d_j} \rho_{j,m} \sum_{p=0}^m \dbinom{m}{p} j^{m-p} \int_{0}^z x^{\beta+n+j-1} \theta^p F(x) \mathrm{d}x.$$ Successive integrations by parts give \begin{multline} \label{eqlemmerec1} \int_{0}^z x^{\beta+n+j-1} \theta^p F(x) \mathrm{d}x = z^{\beta+n+j} \sum_{q=0}^{p-1} (-1)^{p-q-1} (\beta+n+j)^{p-q-1} \theta^q F(z) \\+(-1)^p (\beta+n+j)^p z^{\beta} F_{\beta, n+j}^{[1]}(z).\end{multline} Therefore, diving both sides of the equality by $z^{\beta}$ and using the equality $$\sum\limits_{m=0}^{d_j} \rho_{j,m} \sum\limits_{p=0}^m \dbinom{m}{p} (-1)^p (\beta+n+j)^p j^{m-p}=Q_j(-n-\beta),$$ we obtain $$\sum_{j=0}^{\ell} Q_j(-n-\beta) F_{\beta, n+j}^{[1]}(z)=\sum_{j=0}^\ell z^{n+j} B_{j,n,1,\beta}(\theta) F(z)$$ with $$B_{j,n,1,\beta}=\sum\limits_{q=0}^{d_j-1} b_{j,n,1,q,\beta} X^q \; , \quad \quad b_{j,n,1,q,\beta}=\sum\limits_{m=0}^{d_j} \rho_{j,m} \sum\limits_{p=q+1}^m \dbinom{m}{p} j^{m-p} (\beta+n+j)^{p-q-1} (-1)^{p-q-1}.$$ Multiplying both sides of the equality by $\db^{\mu}$, we see that the coefficients of $\db^{\mu} B_{j,n,1,\beta}(X)$ are algebraic integers which are also polynomials in $n$ with coefficients in $\Oal_\K$ of degree at most $d_j-q-1$. This is the desired conclusion.

\item Let $s \in \N^*$. We assume that the result holds for $s$. We saw in the first point that $$\int_{0}^z x^{\beta-1} F_{\beta, n+j}^{[s]}(x) \mathrm{d}x=z^{\beta} F_{\beta, n+j}^{[s+1]}(z).$$ So, by induction hypothesis, \begin{align*}
\sum\limits_{j=0}^\ell Q_{j,\beta}(-n)& F_{\beta, n+j}^{[s+1]}(z) = \gf{1}{z^\beta} \int_{0}^z \sum\limits_{j=0}^\ell Q_{j,\beta}(-n) x^{\beta-1} F_{\beta, n+j}^{[s]}(x) \mathrm{d} x \\
&= \gf{1}{z^\beta} \sum_{j=0}^\ell \sum\limits_{t=1}^{s-1} \gamma_{j,n,t,s,\beta} \int_{0}^z x^{\beta-1} F_{\beta, n+j}^{[t]}(x) \mathrm{d} x+\gf{1}{z^{\beta}} \sum\limits_{j=0}^\ell \int_{0}^{z} x^{\beta+n+j-1} B_{j,n,s,\beta}(\theta) F(x) \mathrm{d} x \\
&= \sum_{j=0}^\ell \sum\limits_{t=1}^{s-1} \gamma_{j,n,t,s,\beta} F_{\beta, n+j}^{[t+1]}(z)+\sum\limits_{j=0}^\ell \sum_{q=0}^{d_j-s} b_{j,n,s,q,\beta} \gf{1}{z^{\beta}} \int_{0}^z x^{\beta+n+j-1} \theta^q F(x) \mathrm{d} x.
\end{align*}
Finally, Equation \eqref{eqlemmerec1} yields $$\sum_{j=0}^\ell Q_j(-n-\beta) F_{\beta, n+j}^{[s+1]}(z)=\sum_{j=0}^\ell \sum_{t=1}^s \gamma_{j,n,t,s+1,\beta} F_{\beta, n+j}^{[t]}(z)+\sum_{j=0}^\ell z^{n+j} B_{j,n,s+1,\beta}(\theta) F(z)$$ where $$\gamma_{j,n,t,s+1,\beta}=\begin{cases} \gamma_{j,n,t-1,s,\beta}\;\;,  &  2 \leqslant t \leqslant s+1 \\ \sum\limits_{i=0}^{\ell} \sum\limits_{q=0}^{d_j-s} (-1)^q (\beta+n+i)^q b_{i,n,s,q,\beta}\;\;, & t=1 \end{cases}$$ \end{itemize} and $$B_{j,n,s+1,\beta}(X)=\sum_{q=0}^{d_j-s} b_{j,n,s,q,\beta} \sum_{p=0}^{q-1} (-1)^{q-p-1} (\beta+n+j)^{q-p-1} X^p \in \Oal_\K[X]$$ has degree at most $d_j-s-1$.\end{dem}

Lemma \ref{lemmerecurrenceFns} implies the following proposition, which is the main result of this subsection. It provides an inhomogeneous recurrence relation satisfied by the sequence of $G$-functions $\left(F_{\beta,n}^{[s]}(z)\right)_{n \in \N^*,\; 0 \leqslant s \leqslant S}$. The important fact in \eqref{eq:recurrenceFns} is that the length of the summations over $j$ does not depend on $n$.

\begin{prop} \label{prop:recurrenceFns}
Let $m \in \N^*$ be such that $m > f-\ell-\beta$ for every exponent $f$ of $L$ at $\infty$ satisfying $f-\beta \in \N$. Then for any $s, n \geqslant 1$,

\textbf{a)} There exist some algebraic numbers $\kappa_{j,t,s,n,\beta} \in \K$ and polynomials $K_{j,s,n,\beta}(z) \in  \K[z]$ of degree at most $n+s(\ell-1)$ such that \begin{equation}\label{eq:recurrenceFns}F_{\beta,n}^{[s]}(z)=\sum_{t=1}^s \sum_{j=1}^{\ell+m-1} \kappa_{j,t,s,n,\beta} F_{\beta,j}^{[t]}(z)+\sum_{j=0}^{\mu-1} K_{j,s,n,\beta}(z)(\theta^j F)(z).\end{equation}

\textbf{b)} There exists a constant $C_1(F,\beta) >0$ such that the numbers $\house{\kappa_{j,t,s,n,\beta}} ( 1 \leqslant j \leqslant \ell+m-1, \; 1 \leqslant t \leqslant s$), and the houses of the coefficients of the polynomials $K_{j,s,n,\beta}(z), 0 \leqslant j \leqslant \mu-1$ are bounded by $H(F,\beta,s,n)$, with $$ \forall n \in \N^*,\;\; \forall 1 \leqslant s \leqslant S, \quad H(F,\beta,s,n)^{1/n} \leqslant C_1(F,\beta)^S.$$

\textbf{c)} Let $D(F,\beta,s,n)$ denote the least common denominator of the algebraic numbers $\kappa_{j,t,s,n',\beta}$ ($1 \leqslant j \leqslant \ell+m-1$, $1 \leqslant t \leqslant s$, $n' \leqslant n$) and of the coefficients of the polynomials $K_{j,s,n',\beta}(z)$ ($0 \leqslant j \leqslant \mu-1$, $n' \leqslant n$). Then there exists a constant $C_2(F,\beta)>0$ such that $$\forall n \in \N^*,\;\; \forall 1 \leqslant s \leqslant S, \quad D(F,\beta,s,n)^{1/n} \leqslant C_2(F,\beta)^S.$$
\end{prop}

The proof of this proposition is, \emph{mutatis mutandis}, the same as the proof of \cite[Proposition 1, p. 16]{FRivoal}. Indeed, Proposition \ref{prop:Gopbasesol} implies that $L_{\beta}=\sum\limits_{j=0}^{\ell} z^j Q_{j, \beta}(\theta+j)$ is a $G$-operator, which enables us to use \cite[Lemma 2, p. 12]{FRivoal} in order to deduce Proposition \ref{prop:recurrenceFns} from Lemma \ref{lemmerecurrenceFns} above.

However, we will present in Section \ref{sec:explicitC} a precise way to compute the constants $C_1(F,\beta)$ and $C_2(F,\beta)$ which was not given in \cite{FRivoal}. In particular, we will see that the constant
$C_1(F,\beta)$ can be chosen independent of $\beta$.

In the next two subsections and in Section \ref{sec:explicitC}, the index $\beta$ relative to $\kappa_{j,t,s,n,\beta}$ and $K_{j,s,n,\beta}$ will be omitted as there is no ambiguity.

\begin{rqu}
Denote by $\mathcal{E}(\beta)$ the set of exponents $f$ of $L$ at $\infty$ such that $f-\beta \in \N$. Then the best possible value for $m$ is $m=\max \big(\{1\} \cup \{f+1-\ell-\beta, f \in \mathcal{E}(\beta) \}\big)=\ell_0(\beta)-\ell+1$ where \begin{equation}\label{eq:defl0beta}\ell_0(\beta):=\max \left(\{ \ell \} \cup \{f-\beta, f \in \mathcal{E}(\beta) \}\right).\end{equation}

Katz's Theorem (see \cite[Theorem 6.1, p. 98]{Dwork}) ensures that the exponents of $L$ at $\infty$ are all rational numbers. Assume that one of them, denoted by $f$, is nonzero, and set, for all $k \in \N$, $\beta_k:=\left(\mathrm{sign}(f)-k \den(f)\right) |f|$. Then we have $f-\beta_k=k \den(f) |f| \in \N$, so that for all $k$, $$\ell_0(\beta_k) \geqslant \mathrm{den}(f) |f| k \xrightarrow[k \rightarrow + \infty]{} +\infty.$$ 
Likewise, if $0$ is the only exponent of $L$ at $\infty$, then $\beta_k=-k-\ell$ satisfies $\ell_0(\beta_k)=k+\ell \rightarrow +\infty$.
\end{rqu}

\subsection{Study of an auxiliary series} \label{subsec:auxiliaryseries}

As in \cite[p. 24]{FRivoal}, we define an auxiliary series $T_{S,r,n}(z)$, which depends on $\beta$ and turns out to be a linear form with polynomial coefficients in the $F_{\beta,u}^{[s]}(z)$ (Proposition \ref{prop:tsrnetcusn}).

For $S \in \N$ and $r \in \N$ such that $r \leqslant S$, let $$T_{S,r,n}(z)=n!^{S-r} \sum_{k=0}^{\infty} \gf{(k-rn+1)_{rn}}{(k+1+\beta)_{n+1}^S} A_k z^{-k}. $$ This series converges for $|z|>R^{-1}$.

The goal of this part is to obtain various estimates on $T_{S,r,n}(z)$ in order to be able to apply a generalization of Nesterenko's linear independence criterion (\cite[Section 3]{FRivoal}). This will provide the lower bound on the dimension of $\Phi_{\alpha,\beta,S}$ in Theorem \ref{th:chapitre2}. 
The control of the size and the denominator of coefficients appearing in the relation between $T_{S,r,n}(z)$ and the $F_{\beta,u}^{[s]}[z)$ (Lemmas \ref{lem:tailleCusn} and \ref{lem:denomCusn}) will play a central role, but the most tedious part in the original paper of Fischler and Rivoal consisted in the use of the saddle point method in order to obtain an asymptotic expansion of $T_{S,r,n}\left(1/\alpha\right)$ as $n \rightarrow +\infty$ for $0 < |\alpha| < R$. Fortunately, we can adapt their proof with only a few minor changes (Lemma \ref{lem:saddlepoint}). 

\bigskip

\begin{prop} \label{prop:tsrnetcusn}
For $n \geqslant \ell_0(\beta)$, there exist some polynomials $C_{u,s,n}(z) \in \K[z]$ and $\tilde{C}_{u,n}(z) \in \K[z]$ of respective degrees at most $n+1$ and $n+1+S(\ell-1)$ such that $$T_{S,r,n}(z)=\sum\limits_{u=1}^{\ell_0(\beta)} \sum_{s=1}^S C_{u,s,n}(z) F_{\beta,u}^{[s]}\left(\gf{1}{z}\right)+\sum_{u=0}^{\mu-1} \widetilde{C}_{u,n}(z)z^{-S(\ell-1)}(\theta^u F)\left(\gf{1}{z}\right).$$\end{prop}

\begin{dem}
Let us write the partial fraction expansion of \begin{equation}\label{eq:elementsimplesTsrn} R_n(X):= n!^{S-r} \gf{X(X-1) \dots (X-rn+1)}{(X+\beta+1)^S \dots (X+\beta+n+1)^S}=\sum\limits_{j=1}^{n+1} \sum\limits_{s=1}^S \gf{c_{j,s,n}}{(X+\beta+j)^s}, \quad c_{j,s,n} \in \Q,\end{equation} so that $$T_{S,r,n}(z)=\sum_{j=1}^{n+1} \sum_{s=1}^S c_{j,s,n} z^j F_{\beta,j}^{[s]}\left(\gf{1}{z} \right).$$ Then \cite[Lemma 4, p. 24]{FRivoal}, Equation \eqref{eq:elementsimplesTsrn} and Proposition \ref{prop:recurrenceFns} altogether yield $$T_{S,r,n}(z)=\sum\limits_{u=1}^{\ell_0(\beta)} \sum_{s=1}^S C_{u,s,n}(z) F_{\beta,u}^{[s]}\left(\gf{1}{z}\right)+\sum_{u=0}^{\mu-1} \tilde{C}_{u,n}(z)z^{-S(\ell-1)}(\theta^u F)\left(\gf{1}{z}\right)\; ,$$ where $$C_{u,s,n}(z)=c_{u,s,n} z^u + \sum_{j=\ell_0(\beta)+1}^{n+1} \sum_{\sigma=s}^S z^j c_{j,\sigma,n} \kappa_{u,s,\sigma,j}\;,$$ and $$\widetilde{C}_{u,n}(z)=\sum_{j=\ell_0(\beta)+1}^{n+1} \sum_{s=1}^S c_{j,s,n} z^{j+S(\ell-1)} K_{u,s,j}\left(\gf{1}{z}\right).$$\end{dem}

We begin by computing an upper bound on the house of the coefficients of the polynomials $C_{u,s,n}(z)$ and $\tilde{C}_{u,n}(z)$ appearing in Proposition \ref{prop:tsrnetcusn}. 

\begin{lem} \label{lem:tailleCusn}
For any $z \in \Qbar$, we have $$\limsup\limits_{n \rightarrow +\infty}  \left(\max\limits_{u,s} \house{C_{u,s,n}(z)}\right)^{1/n} \leqslant C_1(F,\beta)^S r^r 2^{S+r+1} \max(1, \house{z})$$ and $$\limsup\limits_{n \rightarrow +\infty} \left(\max\limits_{u,s} \house{\widetilde{C}_{u,s,n}(z)}\right)^{1/n} \leqslant C_1(F,\beta)^S r^r 2^{S+r+1} \max(1, \house{z}).$$
\end{lem}

\begin{dem}
We are going to draw inspiration from the proof given in \cite[pp. 6--7]{Rivoal2003}. 

Let $n \in \N^*$, $j_0 \in \{ 1, \dots, n+1 \}$ and $s_0 \in \{ 1, \dots, S \}$. The residue theorem yields $$c_{j_0,s_0,n}=\gf{1}{2 i \pi} \int_{|z+\beta+j_0|=1/2} R_n(z) (z+\beta+j_0)^{s_0-1} \mathrm{d} z$$ where $R_n(z)$ has been defined in \eqref{eq:elementsimplesTsrn}. If $|z+\beta+j_0|=\gf{1}{2}$, we have 
\begin{align*}
\left|(z-rn+1)_{rn}\right|&=\prod_{k=0}^{rn-1} |z-rn+1+k|=\prod_{k=0}^{rn-1} \left|z+\beta+j_0-\left(rn-1-k+\beta+j_0\right)\right| \\
& \leqslant  \prod_{k=0}^{rn-1} \left(\gf{1}{2}+rn-(k+1)+|\beta|+j_0\right) \leqslant  \prod_{k=0}^{rn-1} \left(rn-k+|\beta|+j_0\right)=(|\beta|+j_0+1)_{rn} \\
& \leqslant (\widetilde{\beta}+j_0+2)_{rn}=\gf{(\widetilde{\beta}+j_0+rn+1)!}{(\widetilde{\beta}+j_0+1)!}\;,
\end{align*} with $\widetilde{\beta}:=\lfloor |\beta| \rfloor$, where $\lfloor \cdot \rfloor$ denotes the integer part function. Moreover, \begin{align*}
\left|(z+\beta+1)_{n+1}\right| &= \prod_{k=0}^{n} |z+\beta+k+1|=\prod_{k=0}^{n} \left|z+\beta+j_0-(j_0-k-1)\right| \\
& \geqslant \prod_{k=0}^{n} \left| |j_0-k-1|- \gf{1}{2}\right|=\prod_{k=0}^{j_0-3} \left(j_0-k-1-\gf{1}{2} \right) \times \left(\gf{1}{2}\right)^{3} \times \prod_{k=j_0+1}^{n} \left(k+1-j_0-\gf{1}{2} \right) \\
& \geqslant \gf{1}{8} (j_0-2)! (n-j_0)!\;\;.
\end{align*} Therefore, \begin{align*}
|R_n(z)| &\leqslant n!^{S-r} \gf{(\widetilde{\beta}+j_0+rn+1)!}{(\widetilde{\beta}+j_0+1)! (j_0-2)!^S (n-j_0)!^S} 8^S=\dbinom{\widetilde{\beta}+j_0+rn+1}{\widetilde{\beta}+j_0+1} \times \gf{(rn)!}{(j_0-2)!^S (n-j_0)!^S} 8^S n!^{S-r} \\
&= \dbinom{\widetilde{\beta}+j_0+rn+1}{\widetilde{\beta}+j_0+1} \times \gf{(rn)!}{n!^r} \times \dbinom{n-2}{j_0-2}^S \times n^S (n-1)^S \times 8^S.
\end{align*} Hence $$|c_{j_0,s,n}| \leqslant 2^{\widetilde{\beta}+j_0+rn+1} r^{rn} 2^{S(n-2)} (n(n-1))^S 8^S \left(\gf{1}{2}\right)^S \leqslant 2^{\widetilde{\beta}+(r+1)n+1} r^{rn} 2^{S(n-2)} (n(n-1))^S 8^S \left(\gf{1}{2}\right)^S $$ so that, since the last bound is independent of $j_0$, we get $$\limsup_{n \rightarrow + \infty} \left(\max_{1  \leqslant j \leqslant n+1} |c_{j,s,n}|\right)^{1/n} \leqslant r^r 2^{S+r+1}.$$ The desired result follows from this inequality and from point \textbf{b)} of Proposition \ref{prop:recurrenceFns}.
\end{dem}

The following lemma then provides an upper bound on the denominator of the coefficients of $C_{u,s,n}(z)$ and $\tilde{C}_{u,n}(z)$. We recall that $\db$ is the denominator of $\beta$.

\begin{lem} \label{lem:denomCusn}
Let $z \in \K$ and $q \in \N^*$ be such that $qz \in \Oal_{\K}$. Then there exists a sequence $(\Delta_n)_{n \geqslant 1}$ of positive natural integers such that, for any $u,s$: $$\Delta_n C_{u,s,n}(z) \in \Oal_{\K}, \quad\quad \Delta_n \widetilde{C}_{u,n}(z) \in \Oal_{\K}, \quad \text{and} \quad \lim\limits_{n \rightarrow + \infty} \Delta_n^{1/n}=q C_2(F, \beta)^S \db^{2r} e^S.$$
\end{lem}

\begin{dem}[of  Lemma \ref{lem:denomCusn}]
We are going to follow the proof given by Rivoal in \cite[pp. 7--8]{Rivoal2003}. For practical reasons, we will work with $$\widetilde{R}_n(t)=R_n(t-1)=n!^{S-r} \gf{(t-rn)_{rn}}{(t+\beta)_{n+1}^S}\; ,$$ rather than with $R_n(t)$. For any $j_0 \in \{0, \dots, n \}$, we have \begin{equation} \label{cjsnderivRn} \forall 1 \leqslant s \leqslant S\;, \quad c_{j_0+1, s,n}=D_{S-s}\left(\widetilde{R}_n(t)(t+\beta+j_0)^S\right)_{|t=-\beta-j_0},\end{equation} with $D_{\lambda}=\gf{1}{\lambda!} \gf{\mathrm{d}^{\lambda}}{\mathrm{d}t^{\lambda}}$. Consider the following decomposition: $$\widetilde{R}_n(t)(t+\beta+j_0)^S=\left(\prod_{\ell=1}^r F_{\ell}(t)\right) H(t)^{S-r}$$ with, for $1 \leqslant \ell \leqslant r$, $$F_{\ell}(t)=\gf{(t-n\ell)_n}{(t+\beta)_{n+1}}(t+\beta+j_0)\; , \quad \text{and} \quad H(t)=\gf{n!}{(t+\beta)_{n+1}} (t+\beta+j_0).$$ We obtain $$F_{\ell}(t)=1+\sum_{p=0 \atop p \neq j_0}^n \gf{j_0-p}{t+\beta+p} f_{p,\ell,n} \; , \quad f_{p,\ell,n}=(-1)^{n-p} \dbinom{n}{p} \dbinom{\beta+p+\ell n}{n}$$ and $$H(t)=\sum_{p=0 \atop p \neq j_0}^n \gf{(j_0-p) h_{p,n}}{t+\beta+p}, \quad h_{p,n}=(-1)^p \dbinom{n}{p}.$$ Note that $h_{p,n} \in \N^*$. Hence, if $\lambda \in \N$, $$D_{\lambda}(F_{\ell}(t))_{|t=-\beta-j_0}=\delta_{0,\lambda}+\sum_{p=0 \atop p \neq j_0}^n (-1)^{\lambda} \gf{(j_0-p)f_{p,\ell,n}}{(p-j_0)^{\lambda+1}}=\delta_{0,\lambda}-\sum_{p=0 \atop p \neq j_0}^n  \gf{f_{p,\ell,n}}{(j_0-p)^{\lambda}},$$ where $\delta_{0,\lambda}=1$ if $\lambda=0$ and $0$ else, and $$D_{\lambda}(H(t))_{|t=-\beta-j_0}=-\sum\limits_{p=0}^n \gf{h_{p,n}}{(j_0-p)^\lambda}.$$ Thus, for all $1 \leqslant \ell \leqslant r$ and all $\lambda \in \N$, we have $$d_n^{\lambda} \Delta^{(1)}_n D_{\lambda}(F_{\ell}(t))_{|t=-\beta-j_0} \in \Z \quad \mathrm{and} \quad d_n^{\lambda} D_{\lambda}(H(t))_{|t=-\beta-j_0} \in \Z$$ with $d_n=\mathrm{lcm}(1,2,\dots,n)$ and $\Delta^{(1)}_n \in \N^*$ a common denominator of the $f_{p,\ell,n}$ for all $p,\ell$.

Lemma \ref{lem:denomquotientspochhammer} \textbf{b)} of Subsection \ref{subsec:estimatedenom1/Wbeta} ensures that the integers $\Delta^{(1)}_n=\db^{2n}$ are suitable ones.

Moreover, Leibniz's formula yields $$D_{S-s}\left(\widetilde{R_n}(t)(t+\beta+j_0)^S\right)=\sum_{\boldsymbol{\mu}} D_{\mu_1}(F_1(t)) \dots D_{\mu_r}(F_r(t)) D_{\mu_{r+1}}(H(t)) \dots D_{\mu_{S}}(H(t))\; ,$$ where the sum is on the $\boldsymbol{\mu}=(\mu_1, \dots, \mu_S) \in \N^S$ such that $\mu_1+\dots+\mu_S=S-s$. Finally, using \eqref{cjsnderivRn}, we see that $$ \forall 0 \leqslant j \leqslant n, \; \; \forall 1 \leqslant s \leqslant S \;, \quad d_n^{S-s} \db^{2rn} c_{j+1,s,n} \in \Z.$$ The Prime Number Theorem gives $d_n \leqslant e^{n+o(n)}$ so that the desired conclusion follows from point \textbf{c)} of Proposition \ref{prop:recurrenceFns}. 
\end{dem} 
Let us now explain briefly how the approach of Fischler and Rivoal in \cite{FRivoal} to estimate $T_{S,r,n}\left(1/\alpha\right)$ as $n \rightarrow +\infty$ for $0 < |\alpha| < R$ with the saddle point method can be adapted in our case.

In \cite{FRivoal}, a family of functions $B_1(z), \dots, B_p(z)$ analytic in some half plane $\mathrm{Re}(z) >u$ such that $A(z)=\sum\limits_{j=1}^{p} B_j(z)$ satisfies $A(k)=A_k$ for all large enough integers $k$ has been constructed. Here, $u$ is a positive real number such that $|F(z)|=\mathcal{O}(|z|^{u})$ when $z \rightarrow \infty$ in $\C \setminus (L_0 \cup \dots \cup L_p)$, where the $L_i$ are half-lines (see \cite[p. 28]{FRivoal}). The theory of singular regular points (see \cite[chapter 9]{Hille}) ensures that $u$ exists. Moreover, \cite[Lemma 8, p. 29]{FRivoal} gives, for every $j \in \{1, \dots, p \}$, the following asymptotic expansion of $B_j(tn)$, when $n$ tends to infinity: $$B_j(tn)=\kappa_j \gf{\log(n)^{s_j}}{(tn)^{b_j} \xi_j^{tn}}\left(1+\mathcal{O}\left(\gf{1}{\log(n)}\right)\right),$$ where $s_j \in \N$, $b_j \in \Q$, $\kappa_j \in \C^*$, and $\xi_1, \dots, \xi_p$ are the finite singularities of $F(z)$.  Furthermore, the implicit constant is uniform in any half-plane $\mathrm{Re}(t) \geqslant d, d>0$.

We define $$\mathcal{B}_{S,r,n,j}(\alpha)=\displaystyle\int_{c-i\infty}^{c+i\infty} B_j(tn) \gf{n!^{S-r} \Gamma((r-t)n) \Gamma(tn+\beta+1)^S \Gamma(tn+1)}{\Gamma((t+1)n+\beta+2)^S} (-\alpha)^{tn} \mathrm{d}t\;,$$ for $1 \leqslant j \leqslant p$, where $0<c<r$.

Adapting the computations done in \cite[p. 31]{FRivoal}, based on the residue formula, we obtain the following result:

\begin{lem}
If $0 < |\alpha| < R$ and $r>u$ then for $n$ large enough, we have $$T_{S,r,n}\left(\gf{1}{\alpha}\right)=\sum_{j=1}^p \gf{(-1)^{rn} n}{2 i \pi} \mathcal{B}_{S,r,n,j}(\alpha).$$
\end{lem}

It is now a matter of studying the asymptotic behavior of $\mathcal{B}_{S,r,n,j}(\alpha)$ when $n$ tends to infinity; this is a sensitive step using the saddle point method.

Stirling's formula provides the following asymptotic expansion of $\mathcal{B}_{S,r,n,j}(\alpha)$:
$$\mathcal{B}_{S,r,n,j}(\alpha)=(2\pi)^{(S-r+2)/2} \kappa_j \gf{\log(n)^{s_j}}{n^{(S+r)/2+b_j}} \displaystyle\int_{c-i\infty}^{c+i \infty} g_{j,\beta}(t) e^{n\varphi(-\alpha/\xi_j,t)}\left(1+\mathcal{O}\left(\gf{1}{\log(n)}\right)\right)\mathrm{d}t$$ as $n \rightarrow \infty$, where the constant in $\mathcal{O}$ is uniform in $t$ and  $$g_{j,\beta}(t)=t^{(S+1)/2+S\beta-b_j} (r-t)^{-1/2} (t+1)^{-S(2\beta+3)/2}$$ and $$\varphi(z,t)=t \log(z)+(S+1)t \log(t)+(r-t)\log(r-t)-S(t+1)\log(t+1).$$

Note that $\varphi$ is the same function as in \cite{FRivoal}. Thus, the application of the saddle point method will not change much in this case because $\beta$ appears only in $g_{j,\beta}(t)$. We will have to check that $g_{j,\beta}(t)$ is defined and takes a nonzero value at the saddle point. 

\begin{lem} \label{lem:saddlepoint}
For $z$ such that $0 < |z| <1$ and $-\pi <\arg(z) \leqslant \pi$, let $\tau_{S,r}(z)$ be the unique $t$ such that $\mathrm{Re}(t)>0$ and $\varphi'(z,t)=0$, where $\varphi'(z,t)=\gf{\partial \varphi}{\partial t}(z,t)$.

Assume that $r=r(S)$ is an increasing function of $S$ such that $r=o(S)$ and $Se^{-S/(r+1)}=o(1)$ as $S$ tends to infinity. Then if $S$ is large enough (with respect to the choice of the function $S \mapsto r(S)$), the following asymptotic estimate holds for any $j \in \{1, \dots, p \}$:

$$\mathcal{B}_{S,r,n,j}(\alpha) = (2 \pi)^{(S-r+3)/2} \gf{\kappa_j \gamma_{j,\beta}}{\sqrt{-\psi_j}} \gf{\log(n)^{s_j} e^{\varphi_j n}}{n^{(S+r+1)/2+\beta_j}}(1+o(1)), \quad n \rightarrow +\infty,$$ where $\tau_j=\tau_{S,r}(-\alpha/\xi_j)$, $\varphi_j=\varphi(-\alpha/\xi_j,\tau_j)$, $\psi_j=\varphi''(-\alpha/\xi_j,\tau_j)$ $\gamma_{j,\beta}=g_{j,\beta}(\tau_j)$. Moreover, for any $j \in \{1, \dots, p\}$, we have $\kappa_j \gamma_{j,\beta} \psi_j \neq 0$.
\end{lem}

Note that $\varphi_j, \psi_j, \tau_j$ are the same quantities as in \cite{FRivoal}. The condition on $r$ is in particular satisfied by $r=\left\lfloor S/(\log(S))^2 \right\rfloor$.

\begin{dem}
Only the fourth step of the proof in \cite{FRivoal} has to be adapted to this case in order to apply the saddle point method.

We have $g_{j,\beta}(t)=t^{(S+1)/2+S\beta-b_j} (r-t)^{-1/2} (t+1)^{-S(2\beta+3)/2}$. Hence, denoting $\tau=\tau_{S,r}(z)$, $$g_{j,\beta}(\tau)=\gf{\tau^{(S+1)/2+S\beta}}{\tau^{b_j}(r-\tau)^{1/2}(\tau+1)^{S(2\beta+3)/2}}.$$ But as mentioned in \cite[Step 1, p. 33]{FRivoal}, we have $z \tau^{S+1}-(r-\tau)(\tau+1)^S=0$, so that $r-\tau=\gf{z \tau^{S+1}}{(\tau+1)^S}$, hence $$g_{j,\beta}(\tau)=\gf{\tau^{S\beta-(S+1)/2}}{\tau^{b_j}z^{1/2}(\tau+1)^{S(\beta+1)}} \neq 0$$ because $\mathrm{Re}(z)>0$. \end{dem}

We then deduce from the result above the following proposition, which is the adaptation of \cite[Lemma 7, p. 26]{FRivoal}. The key point, that is proved in \cite[p. 41]{FRivoal}, is that the numbers $e^{\varphi_j}$ are pairwise distinct if we make the additional assumption that $r^{\omega} e^{-S/(r+1)}=o(1)$ for any $\omega >0$. It is satisfied by $r=\left\lfloor S/(\log(S))^2 \right\rfloor$.

\begin{prop} \label{prop:estimationtsrn} Let $\alpha \in \C$ be such that $0< | \alpha| < R$. Assume that $S$ is sufficiently large (with respect to $F$ and $\alpha$), and that $r$ is the integer part of $S/(\log S)^2$. Then there exist some integers $Q \geqslant 1$ and $\lambda\geqslant 0$, real numbers $a$ and $\kappa$, nonzero complex numbers $c_{1,\beta}$,\ldots, $c_{Q,\beta}$, and pairwise distinct complex numbers $\zeta_1$, \ldots, $\zeta_Q$, such that $|\zeta_q|=1$ for any $q$ and
$$
T_{S,r,n}(1/\alpha ) = a^n n^\kappa \log (n)^\lambda \left( \sum_{q=1}^Q c_{q, \beta} \zeta_q^n +o(1)\right) \mbox{ as } n\rightarrow\infty,
$$
and
$$0 < a  \leqslant \gf{1}{r^{S-r}}.$$
\end{prop}

The proof of this result is, \emph{mutatis mutandis} the same as in \cite[pp. 41--42]{FRivoal} but we give a sketch of it for the reader's convenience.

\begin{dem}[(sketch)]
By Lemma 5, we have $$T_{S,r,n}\left(\gf{1}{\alpha}\right)=\sum_{j=1}^p \gf{(-1)^{rn} n}{2 i \pi} \mathcal{B}_{S,r,n,j}(\alpha).$$

\item The asymptotic expansion of $\mathcal{B}_{S,r,n,j}(\alpha)$ provided by Lemma 6 then implies that $$T_{S,r,n}\left(\gf{1}{\alpha}\right) =  \gf{(-1)^{rn}}{2 i \pi} \gf{(2 \pi)^{(S-r+3)/2}}{n^{(S+r-1)/2}} \sum_{j=1}^{p} \gf{\kappa_j \gamma_{j,\beta}}{\sqrt{-\psi_j}} n^{-\beta_j} \log(n)^{s_j} e^{\varphi_j n}(1+o(1)), \quad n \rightarrow +\infty.$$

We consider $J=\{ j_1, \dots, j_Q\}$ the set of the $j \in \{1, \dots, p\}$ such that $(\mathrm{Re}(\varphi_j), -\beta_j-(S+r-1)/2, s_j)$ is maximal for the lexicographic order, equal to some $(a,\kappa,\lambda) \in \R^3$. Then we may neglect the other terms of the sum; precisely, we have 
$$T_{S,r,n}\left(\gf{1}{\alpha}\right)=a^n n^{\kappa} \log(n)^s \sum_{q=1}^{Q} c_{q,\beta}  \zeta_q^n (1+o(1))$$ where $\zeta_q :=\exp\left(i \mathrm{Im}(\varphi_{j_q})\right)$ and $c_{q,\beta} := \gf{(-1)^{rn}}{2 i \pi} (2 \pi)^{(S-r+3)/2} \gf{\kappa_{j_q} \gamma_{j_q,\beta}}{\sqrt{-\psi_{j_q}}}$

Finally, the difficult point is to prove that the $\zeta_q$ are pairwise distinct. This comes from the fact, proved in \cite[p. 42]{FRivoal} that the $\varphi_j$ are pairwise distinct. As we mentioned it above, it is crucial for this purpose to make the assumption that $r^{\omega} e^{-S/(r+1)}=o(1)$ for any $\omega >0$. 
\end{dem}

\subsection{Proof of Theorem  \cal{\ref{th:chapitre2}}} \label{subsec:preuvetheoremeprincipal}

We are now going to prove the main theorem of this paper. The upper bound on the dimension of $\Phi_{\alpha,\beta,S}$ arise from the recurrence relation \eqref{eq:recurrenceFns} of Proposition \ref{prop:recurrenceFns} above. On the other hand, by the estimates of Subsection \ref{subsec:auxiliaryseries}, we can now apply a linear independence criterion à la Nesterenko (\cite[Theorem 4, p. 8]{FRivoal}) to obtain a lower bound.

For the sake of clarity, we reproduce here, with some adaptations, the proof of \cite[pp. 26-27]{FRivoal}. 

Let $\alpha$ be a nonzero element of $\K$ such that $|\alpha| < R$; choose $q \in \N^*$ such that $\gf{q}{\alpha} \in \Oal_{\K}$.  By Lemmas \ref{lem:tailleCusn} and \ref{lem:denomCusn}, $p_{u,s,n}:=\Delta_n C_{u,s,n}(1/\alpha)$ and 
$\tilde{p}_{u,n}:=\Delta_n \widetilde{C}_{u,n}(1/\alpha)$  belong to $\Oal_{\K}$ and for any $u$,~$s$, $$\limsup_{n\rightarrow +\infty} \max_{u,s}(\house{p_{u,s,n}}^{1/n}, \house{\tilde{p}_{u,n}}^{1/n}) \leqslant b:=q  C_1(F,\beta)^S C_2(F,\beta)^S \db^{2r} e^S r^r 2^{S+r+1}\max(1, \house{1/\alpha}).$$
Using Proposition \ref{prop:tsrnetcusn}, we consider 
$$\tau_n:= \Delta_n T_{S,r,n}\left(\gf{1}{\alpha}\right)=\sum\limits_{u=1}^{\ell_0(\beta)} \sum_{s=1}^S p_{u,s,n} F_{\beta,u}^{[s]}(\alpha)+\sum_{u=0}^{\mu-1} \tilde{p}_{u,n} \alpha^{S(\ell-1)}(\theta^u F)(\alpha).$$

Choosing $r = \left\lfloor S/\log(S)^{2}\right\rfloor$,  Lemma \ref{lem:denomCusn} and Proposition \ref{prop:estimationtsrn} yield as $n$ tends to infinity:
  $$ \tau_n = a_0^{n(1+o(1))}  \left( \sum_{q=1}^Q c_{q,\beta} \zeta_q^n +o(1)\right) \mbox{ with } 0 < a_0  < \gf{q C_2(F,\beta)^S \db^{2r} e^S}{r^{S-r}}.$$

Let $\Psi_{\alpha,\beta,S}$ denote the $\K$-vector space spanned by the numbers $F_{\beta,u}^{[s]}( \alpha)$ and $(\theta^v F)( \alpha)$, $1\leqslant u \leqslant \ell_0(\beta)$, $1 \leqslant s \leqslant S$, $0 \leqslant v \leqslant \mu-1$.  

By \cite[Corollary 2, p. 9]{FRivoal}, we get 
$$
\dim_{\K} (\Psi_{\alpha,\beta, S})\geqslant \gf{1}{[\K:\Q]} \left(1 - \frac{\log (a_0) }{\log (b)}\right).$$
Now, as $S$ tends to infinity, 
\begin{equation} \label{eq:preuvethprincipal1} \log(b)  = \log(2eC_1(F,\beta)C_2(F,\beta))S+o(S) \mbox{ and } 
\log(a_0) \leqslant  -S\log(S)+ o(S\log S).\end{equation} Indeed, $$\log(b)=S\log(2eC_1(F,\beta)C_2(F,\beta))+r\left(\log(r)+\log(2)+2\log\left(\db\right)\right)+\log\left(2q\max(1, \house{1/\alpha})\right)$$ and we have $r=o(S)$ and $$r\log(r)=\gf{S}{\log(S)^2} \log\left(\gf{S}{\log(S)^2}\right)(1+o(1))=S\left(\gf{1}{\log(S)}-\gf{2 \log\log(S)}{\log(S)^2}\right)(1+o(1))=o(S).$$ On the other hand, $\log(a_0) \leqslant -(S-r) \log r+ S\log(C_2(F,\beta) e)+2r \log\left(\db\right)+\log(q)$ and \begin{align*}-(S-r)\log(r)=-(S-r)\left(\log\left(\gf{S}{\log(S)^2}\right)+o(1)\right)&=-S\log(S)-2S\log\log(S)+o(S \log(S)) \\
&\leqslant -S \log(S)+o(S \log(S)),\end{align*} which proves the second part of \eqref{eq:preuvethprincipal1}, since $r=o(S)=o(S \log(S))$.

Therefore,
\begin{equation}\label{eq:300}
\dim_{\K} (\Psi_{ \alpha,S})\geqslant \frac{1+o(1)}{  [\K:\Q] \log(2e C_1(F,\beta)C_2(F,\beta))}\log(S) \qquad \text{as} \;\; S\to +\infty.
\end{equation}

Point \textbf{a)} of Proposition \ref{prop:recurrenceFns} with $m =\ell_0(\beta)-\ell+1$ and $z= \alpha$ shows that $\Phi_{\alpha,\beta,S}$, which is the $\K$-vector space spanned by the numbers $F_{\beta,u}^{[s]}( \alpha)$ for  $u \geqslant 1$ and $0 \leqslant s \leqslant S$, is a $\K$-subspace of $\Psi_{\alpha,\beta,S}$. In particular, for any $S \geqslant 0$, 
$$
\dim_{\K} (\Phi_{\alpha,\beta,S}) \leqslant \dim_{\K} (\Psi_{\alpha,\beta,S}) \leqslant \ell_0(\beta) S+\mu \;,  
$$
which proves the right-hand side of the inequality in Theorem \ref{th:chapitre2}. On the other hand, we also have  
$$
\dim_{\K} (\Psi_{\alpha, \beta, S}) \leqslant \dim_{\K} (\Phi_{\alpha,\beta,S}) +\mu
$$ 
so that the lower bound \eqref{eq:300} holds as well with  $\Phi_{\alpha,\beta,S}$ instead of  $\Psi_{\alpha, \beta, S}$ because $\mu$ is independent from $S$. This proves the left hand side of the inequality in Theorem \ref{th:chapitre2} with \begin{equation} C(F,\beta):=\log(2e C_1(F,\beta)C_2(F,\beta)).\end{equation} The rest of the paper is devoted to the computation of $C(F,\beta)$.
\section{Quantitative results on the size of a differential system} \label{sec:sizeGop}

In this section, we explain the notion of size of a differential system and give several quantitative results about it. This will be useful to give an explicit expression of the constant $C(F,\beta)$ of Theorem \ref{th:chapitre2} in Section \ref{sec:explicitC}.

$G$-operators constitute a class of differential operators which contain the minimal operators of the $G$-functions. They satisfy a condition of moderate growth on certain denominators, called the \emph{Galochkin condition} (Definition \ref{def:galochkin} below). Following André's notations (\cite[Section IV.4]{Andre}), we can attach to any $G$-operator $L$ a quantity $\sigma(L)$, called \emph{size} of $L$, encoding this condition.

\subsection{Size of a differential system}\label{subsec:sizediffsyst}

Let $\K$ be a number field, $\mu \in \N^*$ and $G \in \mathcal{M}_{\mu}(\K(z))$. In this subsection, we introduce and give some basic properties of the size $\sigma(G)$ of $G$.

For $s \in \N$, we define $G_s$ as the matrix such that, if $y$ is a vector satisfying $y'=Gy$, then $y^{(s)}=G_s y$. In particular, $G_0$ is the identity matrix. The matrices $G_s$ satisfy the recurrence relation $$\forall s \in \N, \quad G_{s+1}=G_sG+G'_s,$$ where $G'_s$ is the derivative of the matrix $G_s$.

\begin{defi}[Galochkin, \cite{Galochkin74}] \label{def:galochkin}
 The system $y'=Gy$ is said to \emph{satisfy the Galochkin condition} if there exists $T(z) \in \K[z]$ such that $T(z)G(z) \in \mathcal{M}_{\mu}(\K[z])$ and \begin{equation}\label{eq:galochkin} \exists C>0 : \; \forall s \in \N, \quad  q_s \leqslant C^{s+1}.\end{equation} where, for $s \in \N$, $q_s \geqslant 1$ is defined as the least common denominator of all coefficients of the entries of the matrices $T(z)^m \gf{G_m(z)}{m!}$, when $m \in \{ 1, \dots, s \}$.
\end{defi}

Actually, if $y'=Gy$ satisfies the Galochkin condition, every polynomial $T(z)$ with coefficients in $\K$ such that $T(z)G(z) \in \mathcal{M}_{\mu}(\K[z])$ satisfies the condition \eqref{eq:galochkin} (see Proposition \ref{prop:lienq'ssigmaG} \textbf{b)} below).

Chudnovsky's Theorem (proved in \cite{Chudnovsky}) states that if $G$ is the companion matrix of the minimal nonzero differential operator $L$ associated to a $G$-function, then the system $y'=Gy$ satisfies the Galochkin condition. That is why we say that $L$ is a \emph{$G$-operator} (see \cite[pp. 717 -- 719]{AndregevreyI} for a review of the properties of $G$-operators). Following \cite[chapter VII]{Dwork}, we are now going to rephrase this condition in $p$-adic terms.

If $\mathfrak{p}$ is a prime ideal of $\Oal_{\K}$, $|\cdot|_{\mathfrak{p}}$ is the $\mathfrak{p}$-adic absolute value defined on $\K$, with the choice of normalisations given in \cite[p. 223]{Dwork}. We recall that the \emph{Gauss absolute value} associated with $|\cdot|_{\mathfrak{p}}$ is the non-archimedean absolute value $$\fonction{|\cdot|_{\mathfrak{p},\mathrm{Gauss}}}{\K(z)}{\R}{\gf{\sum\limits_{i=0}^{N} a_i z^i}{\sum\limits_{j=0}^{M} b_j z^j}}{\gf{\max\limits_{0 \leqslant i \leqslant N} |a_i|_{\mathfrak{p}}}{\max\limits_{0 \leqslant j \leqslant M} |b_j|_{\mathfrak{p}}}.}$$ 

The absolute value $|\cdot|_{\mathfrak{p},\mathrm{Gauss}}$ naturally induces an eponymous norm on $\mathcal{M}_{\mu,\nu}(\K(z))$, defined for all $H=(h_{i,j})_{i,j} \in \mathcal{M}_{\mu,\nu}(\K(z))$ as $|H|_{\mathfrak{p},\mathrm{Gauss}}=\max_{i,j} |h_{i,j}|_{\mathfrak{p},\mathrm{Gauss}}$. It is called the \emph{Gauss norm}. If $\mu=\nu$, $\mathcal{M}_{\mu}(\K(z))$ endowed with $|\cdot|_{\mathfrak{p},\mathrm{Gauss}}$ is a normed algebra. We now use this notation to define the notion of size of a matrix:

\begin{defi}[\cite{Dwork}, p. 227] \label{def:taillematrice}
Let $G \in \mathcal{M}_{\mu}(\K(z))$. The \emph{size} of $G$ is $$\sigma(G):=\limsup\limits_{s \rightarrow + \infty} \gf{1}{s} \sum\limits_{\mathfrak{p} \in \Spec(\Oal_\K)} h(s, \mathfrak{p})$$ where $$ \forall s \in \N, \quad h(s, \mathfrak{p})= \sup_{m \leqslant s} \log^{+} \left\vert \gf{G_m}{m!} \right\vert_{\mathfrak{p}, \mathrm{Gauss}}, \quad \text{with} \quad \log^{+} : x \mapsto \log\left(\max(1,x)\right). $$ The \emph{size} of $Y=\sum\limits_{k=0}^{\infty} Y_k z^k, Y_k  \in \mathcal{M}_{\mu,\nu}(\K)$ is $$\sigma(Y):=\limsup\limits_{s \rightarrow + \infty} \gf{1}{s} \sum\limits_{\mathfrak{p} \in \Spec(\Oal_\K)} \sup_{k \leqslant s} \log^{+} \left\vert Y_k \right\vert_{\mathfrak{p}}. $$
\end{defi}

If $\alpha_1, \dots, \alpha_n \in \K$ are such that $(\alpha_i)=\gf{\mathfrak{a}_i}{\mathfrak{b}_i}$ where $\mathfrak{a}_i$ and $\mathfrak{b}_i$ are coprime ideals of $\Oal_\K$, we define $\den'(\alpha_1, \dots, \alpha_n)$ as the norm $N_{\K/\Q}(\mathfrak{b})$ of the the smallest common multiple $\mathfrak{b}$ of $\mathfrak{b}_1, \dots, \mathfrak{b}_n$ in the sense of the Dedekind rings.

This denominator need not be the smallest common multiple of the $\alpha_i$ in the classical sense, as the example of $\alpha=(1+i)/2$ in $\K=\Q(i)$ shows: it satisfies $\mathrm{den}(\alpha)=2$ and $\mathrm{den}'(\alpha)=4$, since $(1+i)$ is a prime ideal of $\Z[i]=\Oal_\K$. However, we have \begin{equation}\label{eq:denomprime}  \den(\alpha_1, \dots, \alpha_n) \leqslant \den'(\alpha_1, \dots, \alpha_n) \leqslant (\den(\alpha_1, \dots, \alpha_n))^{[\K:\Q]}.\end{equation}

This alternative notion of denominator turns out to be more useful than the usual one in our context because, as proved in \cite[p. 225]{Dwork}, \begin{equation} \label{eq:denomhauteur} \forall \alpha_1, \dots, \alpha_n \in \K, \quad \sum\limits_{\mathfrak{p} \in \Spec(\Oal_\K)}\sup\limits_{1 \leqslant i \leqslant n} \log^{+} |\alpha_i |_{\mathfrak{p}} = \gf{1}{[\K:\Q]} \log\left( \den'(\alpha_1, \dots, \alpha_n)\right).\end{equation}

\begin{rqu}
If $u(z)=\sum\limits_{n=0}^{\infty} u_n z^n$, we have $$\sigma(u)=\gf{1}{[\K:\Q]} \limsup\limits_{s \rightarrow + \infty} \gf{1}{s} \log(\den'(u_0,\dots, u_s)) \geqslant \gf{1}{[\K:\Q]} \limsup\limits_{s \rightarrow + \infty} \gf{1}{s} \log(\den(u_0,\dots, u_s)),$$ so that $$\limsup\limits_{s \rightarrow + \infty} \den(u_0,\dots, u_s)^{1/s} \leqslant e^{[\K:\Q] \sigma(u)}.$$
\end{rqu}

We first prove the following proposition which provides the link between Galochkin's condition and the size $\sigma(G)$.

\begin{prop} \label{prop:lienq'ssigmaG}
Let $T(z) \in \K[z]$ be such that $T(z)G(z) \in \mathcal{M}_{\mu}(\K[z])$. Define, for $s \in \N$, $q_s$ (resp. $q'_s$), the denominator (resp. the $\mathrm{den}'$) of the coefficients of the entries of the matrices $TG$, $T^2 \gf{G_2}{2}$, $\dots$, $T^s \gf{G_s}{s!}$. Then
\begin{enumerateth}
\item For all integer $s$,  $\gf{1}{s} \log(q_s) \leqslant \gf{1}{s} \log(q'_s) \leqslant [\K:\Q] \gf{1}{s} \log(q_s)$. 
\item Set $$h^{-}(T)=\sum\limits_{\mathfrak{p} \in \Spec(\Oal_\K)} \log^{-} |T|_{\mathfrak{p}, \mathrm{Gauss}} \quad \mathrm{and} \quad h^{+}(T)=\sum\limits_{\mathfrak{p} \in \Spec(\Oal_\K)} \log^{+} |T|_{\mathfrak{p}, \mathrm{Gauss}}$$ with $\log^{-} : x \mapsto \log\left(\min(1,x)\right)$. Then $$[\K:\Q] \sigma(G)+h^{-}(T) \leqslant \limsup\limits_{s \rightarrow + \infty} \gf{1}{s} \log(q'_s) \leqslant [\K:\Q] \sigma(G) + h^{+}(T).$$
\item The differential system $y'=Gy$ satisfies the Galochkin condition if and only if $\sigma(G) < + \infty$ (\cite[p. 228]{Dwork}).
\end{enumerateth}
\end{prop}

\begin{dem}
Point \textbf{a)} is a direct consequence of the inequality \eqref{eq:denomprime}. It implies immediately \textbf{c)}.

\textbf{b)} Equation \eqref{eq:denomhauteur} yields $$\log q'_s=[\K:\Q] \sum\limits_{\mathfrak{p} \in \Spec(\Oal_\K)} \sup\limits_{m \leqslant s} \log^{+} \left\vert\gf{T^m G_m}{m!} \right\vert_{\mathfrak{p}, \mathrm{Gauss}}.$$

Moreover, for $s\in \N^*$ and $m \leqslant s$, we have $$\log^{+} \left\vert T^m \gf{G_m}{m!} \right\vert_{\mathfrak{p}, \mathrm{Gauss}} =\log^{+} \left( |T|^m_{\mathfrak{p}, \mathrm{Gauss}} \left\vert\gf{G_m}{m!} \right\vert_{\mathfrak{p}, \mathrm{Gauss}}\right)\leqslant s \log^{+} |T|_{\mathfrak{p}, \mathrm{Gauss}} + \log^{+} \left\vert \gf{G_m}{m!} \right\vert_{\mathfrak{p}, \mathrm{Gauss}}.$$ Hence $$\sup_{m \leqslant s} \log^{+} \left\vert T^m \gf{G_m}{m!} \right\vert_{\mathfrak{p}, \mathrm{Gauss}} \leqslant  s \log^{+} |T|_{\mathfrak{p}, \mathrm{Gauss}}+h(s,\mathfrak{p}),$$ so that $$\sum\limits_{\mathfrak{p} \in \Spec(\Oal_\K)} \sup\limits_{m \leqslant s} \log^{+} \left\vert\gf{T^m G_m}{m!} \right\vert_{\mathfrak{p}, \mathrm{Gauss}} \leqslant s \sum\limits_{\mathfrak{p} \in \Spec(\Oal_\K)} \sup\limits_{m \leqslant s} \log^{+} |T|_{\mathfrak{p}, \mathrm{Gauss}}+h(s,\mathfrak{p}).$$

Symmetrically, by noticing that $\gf{G_m}{m!}=\left(\gf{1}{T}\right)^m \gf{T^m G_m}{m!}$ and $\log^{+}\left\vert\gf{1}{T}\right\vert_{\mathfrak{p},\mathrm{Gauss}}=-\log^{-}|T|_{\mathfrak{p},\mathrm{Gauss}}$, we have $$h(s,\mathfrak{p}) \leqslant -s\sum\limits_{\mathfrak{p} \in \Spec(\Oal_\K)} \sup\limits_{m \leqslant s} \log^{-} |T|_{\mathfrak{p}, \mathrm{Gauss}}+\sum\limits_{\mathfrak{p} \in \Spec(\Oal_\K)} \sup\limits_{m \leqslant s} \log^{+} \left\vert\gf{T^m G_m}{m!}\right\vert_{\mathfrak{p},\mathrm{Gauss}}.$$ We obtain the desired result by dividing by $s$ and taking the superior limit.\end{dem}

\begin{rqu}
A convenient situation happens when $G \in \mathcal{M}_{\mu}(\Q(z))$ and $T(z) \in \Z[z]$ has at least one coefficient equal to $1$. In that case, Proposition \ref{prop:lienq'ssigmaG} summarises into the equality $$\sigma(G)=\limsup\limits_{s \rightarrow + \infty} \gf{1}{s} \log(q_s).$$
\end{rqu}

\medskip

The following technical lemmas show that the size of a differential system or a differential operator is invariant by change of variable $u=z^{-1}$ and by equivalence of differential systems. We begin by studying the effect of the change of variable $u=z^{-1}$ on  differential systems and differential operators.

\begin{lem} \label{lem:chgevarinfty}
\begin{enumerateth}
\item Let $G \in \mathcal{M}_{\mu}(\K(z))$ and $G_{\infty}(u):= -u^{-2} G(u^{-1})$ be the matrix such that $y'=Gy \ssi \widetilde{y}'=G_{\infty}\widetilde{y}$, where $\widetilde{y}(u)=y(u^{-1})$. Then we have
 \begin{equation} \label{eq:formuleGinfty} \forall s \in \N^*, \quad G_{\infty,s}(u)=(-1)^{s} \sum\limits_{k=1}^{s} \gf{c_{s,k}}{u^{s+k}} G_k\left(\gf{1}{u}\right)\end{equation} where \begin{equation} \label{eq:formulecoeffsAtilde} \forall s \in \N^*,\quad  \forall 1 \leqslant k \leqslant s, \quad c_{s,k}=\binom{s-1}{s-k} \gf{s!}{k!} \in \Z.\end{equation}

 \item  Let $$L=P_{\mu}(z)\left(\gf{\mathrm{d}}{\mathrm{d}z}\right)^{\mu}+P_{\mu-1}(z)\left(\gf{\mathrm{d}}{\mathrm{d}z}\right)^{\mu-1}+\dots+P_{0}(z) \in \K(z)\left[\gf{\mathrm{d}}{\mathrm{d}z}\right],$$ we consider $L_{\infty} \in \K(u)\left[\mathrm{d}/\mathrm{d}u\right]$ such that for all $f$,  $L(f(z))=0 \ssi L_{\infty}(g(u))=0$, where $g(u)=f(u^{-1})$. Then the differential systems $y'=(A_L)_{\infty} y$ and $y'=A_{L_{\infty}} y$ are equivalent over $\Qbar(z)$.
 \item We have
 $$L_{\infty}=P_{\mu,\infty}(u)\left(\gf{\mathrm{d}}{\mathrm{d}u}\right)^{\mu}+P_{\mu-1,\infty}(u)\left(\gf{\mathrm{d}}{\mathrm{d}u}\right)^{\mu-1}+\dots+P_{0,\infty}(u),$$ where $$\forall s \in \{1, \dots, \mu\}, \quad P_{s,\infty}(u)=\sum_{s=k}^{\mu} (-1)^s c_{s,k} u^{s+k} P_s\left(\gf{1}{u}\right) \quad\quad \mathrm{and} \quad P_{0,\infty}(u)=P_0\left(\gf{1}{u}\right).$$
 \end{enumerateth}
\end{lem}

Lemma \ref{lem:chgevarinfty} implies that the size is invariant by change of variable $u=z^{-1}$, as stated in the following lemma.

\begin{lem} \label{lem:sigmainfini=sigma}
We keep the same notations as in Lemma \ref{lem:chgevarinfty}.
\begin{enumerateth}
\item For all $G \in \mathcal{M}_n(\K(z))$, we have $\sigma(G_{\infty})=\sigma(G)$.
\item If $G \in \mathcal{M}_n(\K(z))$ and $P \in \GL_n(\K(z))$, let $P[G]:=PGP^{-1}+P'P^{-1}$ be a matrix defining a differential system $y'=P[G]y$ that is equivalent to $y'=Gy$ over $\Qbar(z)$. Then $\sigma(G)=\sigma(P[G])$ (\cite[Lemma 1, p. 71]{Andre}).
\item For all $M \in \Qbar(z)\left[\mathrm{d}/\mathrm{d}z\right]$, we have $\sigma(M_{\infty})=\sigma(M)$.
\end{enumerateth}
\end{lem}

\begin{dem}[of Lemma \ref{lem:chgevarinfty}]

\textbf{a)} We are going to show \eqref{eq:formuleGinfty} by induction on $s$. Moreover, we will prove that $c_{1,1}=1$ and \begin{equation}\label{eq:relationrecurrencecoeffsGtilde}  \forall s \in \N^*, \; \; \forall 1 \leqslant k \leqslant s+1, \quad c_{s+1,k}=\begin{cases}(s+1) c_{s,1} & \mathrm{si} \; \; k=1 \\ c_{s,k-1}+(s+k)c_{s,k} & \mathrm{if} \; \; 2 \leqslant k \leqslant s \\ c_{s,s} & \mathrm{if} \; \; k=s+1 \end{cases}.\end{equation}
In particular, we will have $c_{s,k} \in \Z$ for all $s,k$.

For $s=1$, it is the definition of $G_{\infty}$. Assume that \eqref{eq:formuleGinfty} is satisfied by $s \in \N^*$ ; then, since $G_{s+1}=G_sG+G'_s$,
\begin{align*}
&G_{\infty,s+1}(u)=G_{\infty,s}(u) G_{\infty}(u)+G_{\infty,s}'(u)=-\gf{1}{u^2} G_{\infty,s}(u) G\left(\gf{1}{u} \right)+ G_{\infty,s}'(u) \\
&=(-1)^{s+1} \sum_{k=1}^{s} \gf{c_{s,k}}{u^{s+k+2}} G_k\left(\gf{1}{u}\right) G\left(\gf{1}{u}\right)+(-1)^{s+1} \sum_{k=1}^{s} \gf{c_{s,k}}{u^{s+k+2}} G'_k\left(\gf{1}{u}\right)+(-1)^{s+1} \sum_{k=1}^{s} \gf{(s+k) c_{s,k}}{u^{s+k+1}} G_k\left(\gf{1}{u}\right) \\
&=(-1)^{s+1} \sum_{k=1}^{s} \gf{c_{s,k}}{u^{s+k+2}} \left(G_k\left(\gf{1}{u}\right) G\left(\gf{1}{u}\right)+G'_k\left(\gf{1}{u}\right)\right)+(-1)^{s+1} \sum_{k=1}^{s} \gf{(s+k)c_{s,k}}{u^{s+k+1}} G_k\left(\gf{1}{u}\right) \\
&=(-1)^{s+1} \sum_{k=1}^{s} \gf{c_{s,k}}{u^{s+k+2}} G_{k+1}\left(\gf{1}{u}\right)+(-1)^{s+1} \sum_{k=1}^{s} \gf{(s+k)c_{s,k}}{u^{s+k+1}} G_k\left(\gf{1}{u}\right) =(-1)^{s+1} \sum_{k=1}^{s+1} \gf{c_{s+1,k}}{u^{s+1+k}} G_{k}\left(\gf{1}{u}\right)
\end{align*} where the coefficients $c_{s+1,k}$ satisfy \eqref{eq:relationrecurrencecoeffsGtilde}. This ends the proof of \eqref{eq:formuleGinfty}.

It remains to prove the explicit formula \eqref{eq:formulecoeffsAtilde} for the coefficients $c_{s,k}$.  
Setting $\widetilde{c}_{s,k}=\binom{s-1}{s-k} \gf{s!}{k!}$, it suffices for this purpose to check that the sequences $(c_{s,k})_{s,k}$ and $(\widetilde{c}_{s,k})_{s,k}$ satisfy the same recurrence relation \eqref{eq:relationrecurrencecoeffsGtilde}. Indeed, we have $c_{1,1}=1=\widetilde{c}_{1,1}$.

Let $s \in \N^*$ and $k \in  \{ 1, \dots, s-2 \}$. Then \eqref{eq:formuleGinfty} implies
\begin{align*}
(s+s-k) \widetilde{c}_{s,s-k}+\widetilde{c}_{s,s-(k+1)}&=(2s-k) \binom{s-1}{k} \gf{s!}{(s-k)!}+\binom{s-1}{k+1} \gf{s!}{(s-k-1)!} \\
\overset{\text{Pascal's triangle}}{=}& (2s-k) \binom{s}{k+1} \gf{s!}{(s-k)!}+\binom{s-1}{k+1}\left[ \gf{s!}{(s-k-1)!}-(2s-k)\gf{s!}{(s-k)!}\right] \\
&= (2s-k) \binom{s}{k+1} \gf{s!}{(s-k)!} - \binom{s-1}{k+1} s \gf{s!}{(s-k)!} \\
&= \gf{s!}{(s-k)!}\left[(2s-k) \binom{s}{k+1}-\gf{s!}{(k+1)!(s-k-2)!}\right] \\
&=\gf{s!}{(k+1)!} \gf{s!}{(s-k)!} \left[\gf{2s-k}{(s-k-1)!} -\gf{1}{(s-k-2)!}\right] \\
&= \gf{s!}{(k+1)!} \gf{s!}{(s-k)!} \gf{s+1}{(s-k-1)!} = \binom{s+1}{k+1} \gf{(s+1)!}{(s+1-(k+1))!}=\widetilde{c}_{s+1, s+1-(k+1)}
\end{align*} 
Moreover, $\widetilde{c}_{s+1,1}=(s+1)!=(s+1) \widetilde{c}_{s,1}$ and $\widetilde{c}_{s+1,s+1}=1=\widetilde{c}_{s,s}$.

Therefore, $\left(\widetilde{c}_{s,h}\right)_{s \in \N^* \atop 0 \leqslant h \leqslant s}$ satisfies the recurrence relation \eqref{eq:relationrecurrencecoeffsGtilde}, which proves \eqref{eq:formulecoeffsAtilde}.

\medskip

\textbf{b)} We have for all $f$,  $L(f(z))=0 \ssi L_{\infty}(g(u))=0$ where $g(u)=f(u^{-1})$. Denote $G=A_L$ the companion matrix of $L$. We set, for every $f$ satisfying $L(f(z))=0$, $$y(z) :={}^t \left(f(z), f'(z), \dots, f^{(\mu-1)}(z)\right)\quad \text{and} \quad w(u) :={}^t \left(g(u), g'(u), \dots, g^{(\mu-1)}(u)\right),$$ as well as $\widetilde{y}(u) :=y(u^{-1})$. The solutions of $h'=G_{\infty} h$ (resp.  $h'=A_{L_{\infty}} h$)  are the vectors $\widetilde{y}(u)$ (resp. $w(u)$) when $L(f(z))=0$.

Equation \eqref{eq:formuleGinfty} above yields for all $s \in \N^*$ \begin{equation}\label{eq:GsGinfinis}\widetilde{y}^{(s)}(u)=G_{\infty,s}(u) \widetilde{y}(u)=(-1)^{s} \sum\limits_{k=1}^{s} \gf{c_{s,k}}{u^{s+k}} G_k\left(\gf{1}{u}\right) y\left(\gf{1}{u}\right)=(-1)^{s} \sum\limits_{k=1}^{s} \gf{c_{s,k}}{u^{s+k}} y^{(k)}\left(\gf{1}{u}\right) \end{equation} so that by taking the first component of the vectors on each side of the equality, we obtain  $$\forall s \in \N^*, \quad g^{(s)}(u)=(-1)^{s} \sum\limits_{k=1}^{s} \gf{c_{s,k}}{u^{s+k}} f^{(k)}\left(\gf{1}{u}\right).$$

Hence $w(u)=P(u) \widetilde{y}(u)$ where $P(u)$ is a lower triangular matrix whose diagonal has only nonzero terms. Thus, $P \in \mathrm{GL}_n(\Qbar(u))$   and $\widetilde{y} \mapsto P \widetilde{y}$ is the bijection from the solution set of $h'=G_{\infty} h$ to the solution set of $h'=A_{L_{\infty}}(u) h$. These two differential systems are therefore equivalent on $\Qbar(z)$.

\medskip

\textbf{c)} Since $\left(G_{\infty}\right)_{\infty}=G$, we may exchange the roles of $y$ and $\widetilde{y}$ in \eqref{eq:GsGinfinis} in order to obtain \begin{equation}\label{eq:GsGinfinis2} y^{(s)}(z)=G_{s}(z) y(u)(-1)^{s} \sum\limits_{k=1}^{s} \gf{c_{s,k}}{z^{s+k}} \widetilde{y}^{(k)}\left(\gf{1}{z}\right).\end{equation} Taking likewise the first component of the vectors on each side of the equality \eqref{eq:GsGinfinis2}, we get $$\forall s \in \N^*, \quad f^{(s)}(z)=(-1)^{s} \sum\limits_{k=1}^{s} \gf{c_{s,k}}{z^{s+k}} g^{(k)}\left(\gf{1}{z}\right).$$
Hence \begin{align*}
    L(f(z))=0 &\ssi \sum\limits_{s=1}^{\mu} P_{s}(z) \sum\limits_{k=1}^{s} \gf{c_{s,k}}{z^{s+k}} g^{(k)}\left(\gf{1}{z}\right)+P_0(z)g\left(\gf{1}{z}\right) =0 \\ &\ssi \sum\limits_{k=1}^{\mu} \left(\sum\limits_{s=k}^{\mu} (-1)^s \gf{c_{s,k}}{z^{s+k}} P_s(z)\right) g^{(k)}\left(\gf{1}{z}\right)+P_0(z)g\left(\gf{1}{z}\right)=0 \\ &\ssi P_{\mu,\infty}(u) g^{(\mu)}(u)+\dots+P_{0,\infty}(u) g(u)=0,
\end{align*}
where the $P_{\ell,\infty}(u)$ are as defined in the statement of the lemma.\end{dem}

\begin{dem}[of Lemma \ref{lem:sigmainfini=sigma}]
\textbf{a)} By Lemma \ref{lem:chgevarinfty} \textbf{a)}, we have for $s \in \N^*$, $$G_{\infty,s}(u)=(-1)^{s} \sum\limits_{k=1}^{s} \gf{c_{s,k}}{u^{s+k}} G_k\left(\gf{1}{u}\right), \quad \text{where} \;\;  c_{s,k}=\dbinom{s-1}{s-k} \gf{s!}{k!} \in \N^*, \quad 1 \leqslant k \leqslant s.$$ If $\mathfrak{p} \in \Spec(\Oal_{\K})$, then $\left|\dbinom{s-1}{s-k}\right|_{\mathfrak{p}} \leqslant 1$, so that $$\left|\gf{G_{\infty,s}}{s!}\right|_{\mathfrak{p}, \mathrm{Gauss}} \leqslant \max_{1 \leqslant k \leqslant s} \left(\left|\gf{G_k(u^{-1})}{k!}\right|_{\mathfrak{p},\mathrm{Gauss}}\right)=\max_{1 \leqslant k \leqslant s} \left(\left|\gf{G_k(u)}{k!}\right|_{\mathfrak{p},\mathrm{Gauss}}\right)$$ because if $h(z) \in \K(z)$, then $|h(z^{-1})|_{\mathfrak{p},\mathrm{Gauss}}=|h(z)|_{\mathfrak{p},\mathrm{Gauss}}$.
    Therefore, $$\max_{1 \leqslant k \leqslant s} \left(\left|\gf{G_{\infty,k}(u)}{k!}\right|_{\mathfrak{p},\mathrm{Gauss}}\right) \leqslant \max_{1 \leqslant k \leqslant s} \left(\left|\gf{G_k(u)}{k!}\right|_{\mathfrak{p},\mathrm{Gauss}}\right),$$ whence $\sigma(G_{\infty}) \leqslant \sigma(G)$. Since $\left(G_{\infty}\right)_{\infty}=G$, we finally have $\sigma(G)=\sigma(G_{\infty})$. 
    
\medskip

\textbf{b)} The proof of this statement can be found in \cite[Lemma 1, p. 71]{Andre}.

\medskip

\textbf{c)} If $G=A_{M}$ is the companion matrix of $M$ and $G_{\infty}$ is defined as in the point \textbf{a)}, then $\sigma(G)=\sigma(G_{\infty})$. But Lemma \ref{lem:chgevarinfty} \textbf{b)} implies that the differential systems $y'=A_{M_{\infty}}y$ and $y'=G_{\infty} y$ are equivalent over $\Qbar(u)$, hence, by point \textbf{b)}, $\sigma(M_{\infty})=\sigma(M)$.\end{dem}

\subsection{A quantitative form of a generalization of Chudnovsky's Theorem}\label{subsec:chudnovskyeffectif}

The goal of this subsection is to state and generalize results from the book of Dwork (\cite[chapters VII and VIII]{Dwork}) formulating a relation between the size of $\sigma(G)$ and the size of the $G$-functions that are solutions of the differential system $y'=Gy$. This will be a key point of the method of Section \ref{sec:explicitC}.
Recall that if $\zeta \in \K$, the \emph{absolute logarithmic size} of $\zeta$ is $$h(\zeta)=\sum\limits_{\tau : \K \hookrightarrow \C} |\zeta|_{\tau}+\sum\limits_{\mathfrak{p} \in \Spec(\Oal_{\K})} |\zeta|_{\mathfrak{p}}$$ with, for any embedding $\tau : \K \hookrightarrow \C$, $$ |\zeta|_{\tau}=\begin{cases} |\tau(\zeta)|^{1/[\K:\Q]} \quad \mathrm{if} \; \tau(\K) \subset \R \\ |\tau(\zeta)|^{2/[\K:\Q]} \quad \mathrm{else}. \end{cases}$$

The following theorem is a combination of Theorem 2.1, p. 228, Theorem 3.3, p. 238 and Theorem 4.3, p. 243 in \cite[chapter VII]{Dwork}.

\begin{Th}[\cite{Dwork}] \label{th:liensigmaYsigmaG}
If $G \in \mathcal{M}_{\mu}(\K(z))$ satisfies the \emph{Galochkin condition} $\sigma(G) < + \infty$, then there exists a fundamental matrix of solutions of the system $y'=Gy$ of the form $Y(z) z^C$, where $C \in \mathcal{M}_{\mu}(\Qbar)$ is a matrix with eigenvalues in $\Q$ and $Y(z) \in \mathcal{M}_{\mu}(\K\llbracket z \rrbracket)$. Moreover, $$\sigma(Y) \leqslant \Lambda(G):=\mu^2 \sigma(G)+\mu^2+\mu-1+\mu(\mu-1) H(N_C) + (\mu^2+1) \sum\limits_{\zeta \neq 0, \infty} h(\zeta)$$ where the sum is over the non apparent singularities $\zeta$ of the system $y'=Gy$, $N_C$ is the common denominator of the eigenvalues of $C$ and $$H(N):=\gf{N}{\varphi(N)} \sum\limits_{(j, N)=1} \gf{1}{j}\;,\qquad \text{where} \;\; \varphi \;\; \text{denotes Euler's totient function}.$$ 
\end{Th}

\begin{rqus}
\begin{itemizeth}
\item In the case where $G$ is a companion matrix $A_L$ associated with a differential operator $L \in \K(z)\left[\mathrm{d}/\mathrm{d}z\right]$, $N_C$ is the common denominator of the exponents of $L$ at $0$.
\item For simplicity, we set $\Lambda(L):=\Lambda(A_L)$ in that case.
\end{itemizeth}
\end{rqus}

We may be interested in finding a simple upper bound on the size $\sigma(G)$ appearing in the expression of $\Lambda(G)$ in Theorem \ref{th:liensigmaYsigmaG} above. It is a constant that can be hard to compute, whereas it might be easier to study the behavior of a particular solution of $y'=Gy$. That is the interest of the following quantitative version of Chudnovsky's Theorem.

\begin{Th}[\cite{Dwork}, p. 299] \label{th:chudnovskyeffectif}
Let $y(z) \in \K\llbracket z \rrbracket^{\mu}$ be a vector of $G$-functions and $G \in \mathcal{M}_{\mu}(\K(z))$ a matrix such that $y'=Gy$. Let $T(z) \in \K[z]$ be such that $T(z)G(z) \in \mathcal{M}_{\mu}(\K[z])$, and $t:=1+\max(\deg(TG), \deg(T))$. Then, if $\mu \geqslant 2$, $$\sigma(G) \leqslant (5 \mu^2 t-1-(\mu-1)t)\overline{\sigma}(y)$$ where $\overline{\sigma}(y)=\sigma(y)+\limsup\limits_{s \rightarrow + \infty} \gf{1}{s} \sum\limits_{\tau : \K \hookrightarrow \C} \sup\limits_{m \leqslant s} \log^{+}|y_m|_{\tau}$. If $\mu=1$, we have $$\sigma(G) \leqslant (6 t-1)\overline{\sigma}(y).$$
\end{Th}

\begin{rqus}
\begin{itemizeth}
\item In particular, if $G$ is the companion matrix of a minimal differential operator $L$ of order $\mu \geqslant 2$ of a $G$-function $F(z)$, Theorem \ref{th:chudnovskyeffectif} can be reformulated as $$\sigma(L) \leqslant (5 \mu^2 (\delta+1)-1-(\mu-1)(\delta+1))\overline{\sigma}(F), \quad \delta=\deg_z(L),$$ or if $\mu=1$, $\sigma(L) \leqslant (6(\delta+1)-1) \overline{\sigma}(F)$.

\item In order to estimate the analytic term of $\overline{\sigma}(y)$ in practice, we use the Hadamard-Cauchy formula. Denoting by $R_{\tau}$ the radius of convergence of $\sum\limits_{m \geqslant 0} \tau(y_m) z^m$ for every embedding $\tau : \K \hookrightarrow \C$, we obtain $$\limsup\limits_{s \rightarrow + \infty} \gf{1}{s}  \sum_{\tau : \K \hookrightarrow \C} \max_{m \leqslant s} \log^+ |y_m|_{\tau} \leqslant \gf{1}{[\K:\Q]}\ \sum\limits_{\tau : \K \hookrightarrow \C} \varepsilon_{\tau} \max\big(0,-\log(R_{\tau})\big)$$ with $\varepsilon_{\tau}=1$ if $\tau(\K) \subset \R$ and $2$ else.
\end{itemizeth}
\end{rqus}

\medskip

André proved in \cite{AndregevreyI} an analogue of Chudnovsky's Theorem for Nilsson-Gevrey series of arithmetic type of order $0$. We proved in \cite{LepetitSize} the following theorem, which is a quantitative version of this result, and thus a generalization for Nilsson-Gevrey series of arithmetic type of Theorem \ref{th:chudnovskyeffectif}. It relies on previous work by André (\cite[Section IV.4]{Andre}) on the size of differential modules. The "qualitative" part \textbf{a)} is due to André, and it is "quantified" in \textbf{b)}, which is the main result of \cite{LepetitSize}. Its proof can be found in this paper.

 \begin{Th}\label{th:chudnovskynilssongevreyandre}
 Let $S \subset \Q \times \N \times \N$ a finite set, $(c_{\alpha, k,r})_{(\alpha, k,r) \in S} \in \left(\C^*\right)^{S}$ and $(f_{\alpha, k,r}(z))_{(\alpha, k,r) \in S}$ a family of nonzero $G$-functions. We consider $$f(z)=\sum\limits_{(\alpha, k, r) \in S}^{} c_{\alpha, k,r} z^{\alpha} \log(z)^k f_{\alpha, k,r}(z)$$ a Nilsson-Gevrey series of arithmetic type of order $0$.
 Then:
\begin{enumerateth}
\item The function $f(z)$ is solution of a nonzero linear differential equation with coefficients in $\Qbar(z)$ and the minimal operator $L$ of $f(z)$ over $\Qbar(z)$ is a $G$-operator.(\cite[p. 720]{AndregevreyI})

\item For every $(\alpha,k,r) \in S$, let $L_{\alpha,k,r} \neq 0$ denote a nonzero minimal operator of the $G$-function $f_{\alpha,k,r}(z)$ over $\Qbar(z)$. We set $\kappa$ the maximum of  the integers $k$ such that $(\alpha, k, r) \in S$ for some $(\alpha, r) \in \Q \times \N$ and $A :=\{ \alpha \in \Q : \exists (k, r) \in \N^2, (\alpha, k, r) \in S \}$. Then we have \begin{multline*}
    \quad \sigma(L) \leqslant \max\left(1+\log(\kappa+2), \; 2\big(1+\log(\kappa+2)\big)\log\big(\max\limits_{\alpha \in A} \den(\alpha)\big),\right.\\  \left.\max\limits_{(\alpha, k, r) \in S}\big((1+\log(k+2))\sigma(L_{\alpha, k,r})\big) \right).\quad
\end{multline*}
\end{enumerateth}
\end{Th}

Thus, the combination of Theorem \ref{th:chudnovskynilssongevreyandre} and of Chudnovsky's Theorem (Theorem \ref{th:chudnovskyeffectif}) provides, with the notations of Theorem \ref{th:chudnovskynilssongevreyandre}, a relation between $\sigma(L)$ and the $\overline{\sigma}(f_{\alpha,k,r})$, for $(\alpha,k,r) \in~S$. 

\section{Computation of the constants $C_1$ and $C_2$}\label{sec:explicitC}

The goal of this section is to prove Theorem \ref{th:calculC1etC2} stated in the Introduction. The results of Section  \ref{sec:sizeGop} will be heavily used.

\medskip

To this end, we recall that if $C_1(F,\beta)$ and $C_2(F,\beta)$ satisfy Proposition \ref{prop:recurrenceFns} \textbf{b)} and \textbf{c)}, then \eqref{eq:expressionCFbeta} in Subsection 2.3 implies that $C(F,\beta)=\log(2eC_1(F,\beta) C_2(F,\beta))$ is a suitable constant in Theorem~\ref{th:chapitre2}

In \cite{FRivoal}, the existence of constants $C_1(F,\beta)$ and $C_2(F,\beta)$ for $\beta=0$ of Proposition \ref{prop:recurrenceFns} was proved but no formula nor algorithm was given to compute them explicitly. We want to give an explicit expression of these constants.

Let us first explain how $C_1(F,\beta)$ and $C_2(F,\beta)$ constants are defined. We take $m$ such that for all $n \geqslant m$, $Q_{\ell}(-n-\beta) \neq 0$ and $Q_0(-n-\beta) \neq 0$. A suitable $m$ is the smallest positive integer such that $m > -e-\beta$ and $m>f-\ell-\beta$ for every exponent $e$ of $L$ at $0$ such that $e+\beta \in \Z$ and every exponent $f$ of $L$ at infinity such that $f-\beta \in \N$.

Let $(u_{1,\beta}(n))_{n \geqslant m}, \dots (u_{\ell,{\beta}}(n))_{n \geqslant m}$ be a basis of solution of the homogeneous linear recurrence associated with the operator $L_{\beta, \infty}$ obtained from $L_{\beta}$ by a change a variable $u=1/z$: \begin{equation} \label{eq:recurrenceLbeta}\forall n \geqslant m \;, \quad \sum\limits_{j=0}^{\ell} Q_j(-n-\beta) u(n+j)=0.\end{equation}

\begin{rqu}
If $\ell=0$, the only solution of \eqref{eq:recurrenceLbeta} is the sequence that is null from the index $m$. Hence the only power series $y(z) \in \K\llbracket z \rrbracket$ solutions of $L_{\beta,\infty}(y(z))=0$ are polynomials with coefficients in $\K$.

The converse is true : if $L$ is an operator admitting at least one power series solutions around $0$ and such that the only $y(z) \in \Qbar\llbracket z \rrbracket$ satisfying $L(y(z))=0$ are polynomials, then $\ell=0$.
\end{rqu}

We now assume that $\ell \geqslant 1$. Define $$W_{\beta}(n)=\left \vert \begin{matrix}
u_{1,\beta}(n+\ell-1)&\cdots & u_{\ell,\beta}(n+\ell-1)
\\
u_{1,\beta}(n+\ell-2)&\cdots &u_{\ell,\beta}(n+\ell-2)
\\
\vdots & \vdots & \vdots 
\\
u_{1,\beta}(n)&\cdots &u_{\ell,\beta}(n)
\end{matrix}\right\vert$$ the wronskian determinant associated with the basis $(u_{1,\beta}, \dots, u_{\ell,\beta})$ and \begin{equation}\label{eq:defDjbeta} D_{j,\beta}(n)=(-1)^j\left \vert \begin{matrix}
u_{1,\beta}(n+\ell-2)&\cdots &u_{j-1,\beta}(n+\ell-2)&u_{j+1,\beta}(n+\ell-2) &\cdots & u_{\ell,\beta}(n+\ell-2)
\\
\vdots & \vdots & \vdots & \vdots &\vdots&\vdots 
\\
u_{1,\beta}(n)&\cdots &u_{j-1,\beta}(n)&u_{j+1,\beta}(n)&\cdots & u_{\ell,\beta}(n)
\end{matrix}\right\vert\end{equation} one of its minors of order $j$.

The same argument as in \cite[p. 20]{FRivoal} shows that the constant $C_2(F,\beta)$ (resp. $C_1(F,\beta)$) can be taken as any upper bound on $\limsup\limits_{n \rightarrow + \infty} \delta_n^{3/n}$ (resp. $\limsup\limits_{n \rightarrow + \infty} M_n^{1/n}$), where $\delta_n$ (resp. $M_n$) is a common denominator (resp. the maximum of the absolute value) of the numbers $$\gf{1}{W_{\beta}(k)},\;\; \gf{D_{j,\beta}(k)}{Q_{\ell}(1-k-\beta)}, \;\; u_{j,\beta}(k), \quad k \in \{m, \dots, n \}, \; \; j \in \{1, \dots, \ell \}.$$ The next four subsections are devoted to the computation of suitables constants $C_1(F,\beta)$ and $C_2(F,\beta)$. This will prove Theorem \ref{th:calculC1etC2}. In particular, we will see that $C_1(F,\beta)$ ultimately only depends on $F$ and not on $\beta$.

\bigskip

\textit{Preliminary remark.}
The method we are going to present to compute $C_2(F,\beta)$ can be refined when $\ell=1$. In this case, there exists a much simpler and more direct way of proceeding.

Indeed, if $\ell=1$, then \eqref{eq:recurrenceLbeta} implies that $$\forall n \geqslant m, u(n+1)=\gf{Q_0(-n-\beta)}{Q_1(-n-\beta)} u(n)$$ so that, denoting $Q_0(X)=\gamma_0 \prod\limits_{i=1}^{\mu} (X-e_i)$ and $Q_{\ell}(X)=\gamma_{\ell} \prod\limits_{i=1}^{\mu} (X+f_i-\ell)$, we obtain $$u(n)=\left( (-1)^{\ell} \gf{\gamma_{\ell}}{\gamma_0} \right)^{n-m} \prod_{i=1}^{\mu} \gf{\left(m-f_i+\ell+\beta\right)_{n-m}}{\left(m+e_i-\beta\right)_{n-m}}.$$ The same computation as in Subsection \ref{subsec:estimatedenom1/Wbeta} below then shows that 
\begin{multline}\label{eq:expressionC2ell=1}\limsup \delta_n^{3/n} \leqslant \den(1/\gamma_0)^3 \den(1/\gamma_{\ell})^3 \den(\mathbf{e},\beta)^{6 \mu} \den(\mathbf{f},\beta)^{6 \mu} \\ \exp\left(3(\mu+1) \den(\mathbf{f},\beta) + 3\mu \den(\mathbf{e},\beta)\right):= C_2(F,\beta).\end{multline} The constant $C_2(F,\beta)$ thus defined is a suitable one in Proposition \ref{prop:recurrenceFns} in the case $\ell=1$.

In everything that follows, we now assume that $\ell  \geqslant 2$. Note that Theorem \ref{th:calculC1etC2} is still true (but of lesser interest) if $\ell=1$.

\subsection{Estimate of the denominator of the $u_{j,\beta}(n)$}\label{subsec:estimatedenomujbeta}

In this subsection, we will rely on results of Section \ref{sec:sizeGop} about the size of the $G$-operators in order to estimate the denominator of $u_{j,\beta}(n)$ for $j \in \{1, \dots, \ell \}$ as $n$ tends to infinity. Here is how we will proceed.

We can show as in \cite[p. 13]{FRivoal}, that if $(u_n)_{n \geqslant m}$ is a solution of \eqref{eq:recurrenceLbeta} (in particular, $(u_{j,\beta}(n))_{n \geqslant m}$ is in that case), then $U(z^{-1})=\sum\limits_{n=m}^{\infty} u_n z^{-n}$ is a solution of $\widetilde{L}_{\beta}(U(z))=0$, with \begin{equation} \label{eq:defLbetatilde} \widetilde{L}_{\beta}=\left(\ddz\right)^{\ell} z^{m-1} L_{\beta}\; ,\end{equation} \emph{i.e.} $U(z)$ is a solution of the operator $\widetilde{L}_{\beta,\infty}$ obtained by the change of variable $u=z^{-1}$.

The goal of this subsection is to prove the following proposition and to compute the constant $\Lambda_0(L,\beta)$ appearing in it in terms of the parameter $\beta$, the function $F$ and its minimal operator $L$. 

\begin{prop}\label{prop:denomujbeta}
Denote, for $0 \leqslant j \leqslant \ell, Q_j(X)=\sum\limits_{k=0}^{\mu} q_{j,k} X^k$ . Then we have
\begin{equation}\label{eq:denomsection2.1}\limsup\limits_{n \rightarrow + \infty} \den\big(u_{j,\beta}(k), \; 1 \leqslant j \leqslant \ell, \; m \leqslant k \leqslant n\big)^{1/n} \leqslant \exp\left([\K:\Q] \Lambda_0(L,\beta)\right).\end{equation} where 
\begin{multline}\Lambda_0(L,\beta)=(\mu+\ell)^2 \left( \ell + 1+\log(2)\right) \max\left(1, 2 \log(\db), \sigma(L) \right)+(\mu+\ell)^2+\mu+\ell-1\\+(\mu+\ell)(\mu+\ell-1) H\left(\den(f_1,\dots,f_{\ell}) \db\right)+ \left((\mu+\ell)^2+1\right) \sum\limits_{\zeta \neq 0, \infty} h\left(1/\zeta \right)\end{multline} where the last sum is on the roots $\zeta$ of $\chi_L(z):=\sum\limits_{j=0}^{\ell} q_{j,\mu} z^j$.
\end{prop}

The following simple lemma will be useful for the proof of Proposition \ref{prop:denomujbeta}, since it allows us to find the singularities of $\widetilde{L}_{\beta,\infty}$ in terms of the singularities of $L_{\beta}$. We give a proof for the reader's convenience. 

\begin{lem}  \label{lem:singulariteschgtvarinfini}
Let $M \in \Qbar(z)\left[\mathrm{d}/\mathrm{d}z\right]$. Then the singularities $\xi \in \C^*$ of $M_{\infty}$ are the $\zeta^{-1}$, when $\zeta$ is a finite nonzero singularity of $M$.
\end{lem}
\begin{dem}[of Lemma \ref{lem:singulariteschgtvarinfini}]
We set $$M=P_{\mu}(z)\left(\gf{\mathrm{d}}{\mathrm{d}z}\right)^{\mu}+P_{\mu-1}(z)\left(\gf{\mathrm{d}}{\mathrm{d}z}\right)^{\mu-1}+\dots+P_{0}(z).$$ Then, by Lemma \ref{lem:chgevarinfty} \textbf{c)}, we have  $$M_{\infty}=P_{\mu,\infty}(u)\left(\gf{\mathrm{d}}{\mathrm{d}u}\right)^{\mu}+P_{\mu-1,\infty}(u)\left(\gf{\mathrm{d}}{\mathrm{d}u}\right)^{\mu-1}+\dots+P_{0,\infty}(u),$$ where \begin{equation}\label{eq:singulariteschgtvarinfini-eq1}\forall \ell \geqslant 1, \quad P_{\ell,\infty}(u)=\sum_{k=\ell}^{\mu} (-1)^k \dbinom{k-1}{k-\ell} \gf{k!}{\ell!} u^{\ell+k} P_k\left(\gf{1}{u}\right) \quad\quad \mathrm{and} \quad P_{0,\infty}(u)=P_0\left(\gf{1}{u}\right).\end{equation}

The singularities of $M_{\infty}$ are the $\xi$'s such that there exists $\ell \in \{0, \dots, \mu-1\}$ such that $P_{\ell,\infty}/P_{\mu,\infty}$ has a pole at $\xi$.

We write for all $k$ \begin{equation} \label{eq:singulariteschgtvarinfini-eq2}
P_k(z)=\alpha_k u^{e_k} \prod\limits_{i=1}^{d_k} (z-r_{i,k})^{\tau_{i,k}}.\end{equation} If $\zeta \neq 0, \infty$ is not a singularity of $M$, then for all $k$, $\zeta$ is not a pole of $P_k/P_{\mu}$. In particular, even in the case when $\zeta$ is a zero of $P_{\mu}(z)$, it is a zero of higher order of $P_k$ for all $k \in \{0, \dots, \mu-1\}$

Observing that if $r \in \C^*$, $$\gf{1}{u}-r=-\gf{1}{ru}\left(u-\gf{1}{r}\right),$$ the expression \eqref{eq:singulariteschgtvarinfini-eq2} implies the existence of $\gamma_k \in \C$ and $f_k \in  \Z$ such that $$\gf{P_k(1/u)}{P_{\mu}(1/u)}=\gamma_k u^{f_k} Q_k(u),$$ where $\zeta^{-1}$ is not a pole of $Q_k(u) \in \Qbar(u)$. Finally, since $P_{\mu,\infty}(u)=(-1)^{\mu} u^{2\mu} P_{\mu}(u^{-1})$, Equation \eqref{eq:singulariteschgtvarinfini-eq1} shows that for all $\ell \in \{0, \dots, \mu-1 \}$, $P_{\ell,\infty}(u)/P_{\mu,\infty}(u)$ doesn't have a pole at $\zeta^{-1}$. In other words, $\zeta^{-1}$ is an ordinary point of $M_{\infty}$.

Since $(M_{\infty})_{\infty}$ is equal to $M$ up to multiplication by an element of $\Qbar(z)$, we obtain by symmetry that the ordinary points $\xi \neq 0, \infty$ of $M_{\infty}$ are the $\zeta^{-1}$, when $\zeta \neq 0, \infty$ is an ordinary point of $M$. Therefore, by definition, the singularities of $M_{\infty}$ are the $\zeta^{-1}$, when $\zeta \neq 0, \infty$ is a singularity of $M$.
\end{dem}
\begin{dem}[of Proposition \ref{prop:denomujbeta}]

As mentioned at the beginning of this subsection, if $1 \leqslant j \leqslant \ell$, then $U(z)=\sum\limits_{n=m}^{\infty} u_{j,\beta}(n) z^{n}$ is a solution of $\widetilde{L}_{\beta,\infty}(y(z))=0$, where $\widetilde{L}_{\beta,\infty}$ is defined by Equation \eqref{eq:defLbetatilde}. Thus, by Theorem \ref{th:liensigmaYsigmaG} and the remark after Equation \eqref{eq:denomhauteur}, we have \begin{equation}\limsup\limits_{n \rightarrow + \infty} \den\left(u_{j,\beta}(k), \; 1 \leqslant j \leqslant \ell, \; m \leqslant k \leqslant n\right)^{1/n} \leqslant \exp\left([\K:\Q] \Lambda\left(\widetilde{L}_{\beta,\infty}\right)\right).\end{equation}

The proof consists now essentially in bounding the constant $\Lambda(\widetilde{L}_{\beta,\infty})$ of Theorem \ref{th:liensigmaYsigmaG} in terms of $\beta$, $L$ and $F$.

\begin{itemize}
    \item Since $L_{\beta}$ has order $\mu$, the order of $\widetilde{L}_{\beta,\infty}$ is $\mu+\ell$.
    \item The set of the exponents of $\widetilde{L}_{\beta,\infty}$ at $0$ is the set of the exponents of $\widetilde{L}_{\beta}$ at $\infty$, which is the union of the set of the exponents of $L_{\beta}$ at $\infty$ and of the set of the exponents of $\left(\mathrm{d}/\mathrm{d}z\right)^{\ell} z^{m-1}$ at $\infty$, that are all integers. Therefore the constant $N_C$ in Theorem \ref{th:liensigmaYsigmaG} is in that case $\den(f_1+\beta, \dots, f_{\ell}+\beta)$, which divides $\den(f_1,\dots, f_{\ell}) \db$.
    
\item Lemma \ref{lem:singulariteschgtvarinfini} ensures that the singular points $\xi \neq 0, \infty$ of $\widetilde{L}_{\beta,\infty}$ are exactly the $\zeta^{-1}$, when $\zeta$ is a finite nonzero singularity of $\widetilde{L}_{\beta}$. 
    
    But $\widetilde{L}_{\beta}=\left(\mathrm{d}/\mathrm{d}z\right)^{\ell} z^{m-1} L_{\beta}$, so Leibniz's formula shows that if $$L_{\beta}=\sum\limits_{k=0}^{\mu} P_{k,\beta}(z) \left(\ddz\right)^{k}\;, \quad \mathrm{with} \quad P_{k,\beta}(z) \in \K[z],$$ then $\widetilde{L}_{\beta}=z^{m-1}P_{\mu,\beta}(z) \left(\mathrm{d}/\mathrm{d}z\right)^{\mu+\ell}+M$, where $M$ is a differential operator of order less than $\mu+\ell$. Hence the nonzero finite singularities of $\widetilde{L}_{\beta}$ are among the roots of $P_{\mu,\beta}(z)$.
    
    Moreover, we can write \begin{align*}\alpha z^{\mu-\omega} L &=\sum_{j=0}^{\ell} z^j Q_j(\theta+j)=\sum_{j=0}^{\ell} z^j \sum_{k=0}^{\mu} q_{j,k} (\theta+j)^k \\ &=\sum_{j=0}^{\ell} z^j \sum_{k=0}^{\mu} q_{j,k} \sum_{s=0}^k \dbinom{k}{s} j^{k-s} \theta^s=\sum_{s=0}^{\mu} \left(\sum_{j=0}^{\ell} z^j \sum_{k=s}^{\mu} q_{j,k} j^{k-s}\right) \theta^s \end{align*} and likewise $$\gf{1}{\db^{\mu}} L_{\beta}=\sum_{j=0}^{\ell} z^j Q_j(\theta+j-\beta)=\sum_{s=0}^{\mu} \left(\sum_{j=0}^\ell z^j \sum_{k=s}^{\mu} q_{j,k} (j-\beta)^{k-s}\right) \theta^s.$$
    
    The coefficient in front of $\theta^{\mu}$ happens to be the same for $L$ and $L_{\beta}$; it is equal to $\chi_L(z):=\sum\limits_{j=0}^{\ell} q_{j,\mu} z^j$ and it is independent of $\beta$.
    Therefore, since $$\theta^{\mu}=z^{\mu} \left(\ddz\right)^{\mu}+\sum\limits_{k=1}^{\mu-1} a_k z^k \left(\ddz\right)^k \; , \quad a_k \in \C,$$ we have $L_{\beta}=z^{\mu} \chi_L(z) (\mathrm{d}/\mathrm{d}z)^{\mu}+M_{\beta}$, where $M_{\beta}$ has order less than $\mu$. 
    
    So finally, $P_{\mu,\beta}(z)=z^{\mu} \chi_L(z)$ and the nonzero finite singularities of $\widetilde{L}_{\beta}$ are among the roots of $\chi_L(z)$ which does not depend on $\beta$. It may happen that the set of singularities of $\widetilde{L}_{\beta}$ varies with $\beta$ inside the set of roots of $\chi_L$, but this is not important for us.

    \item We have $\sigma(\widetilde{L}_{\beta,\infty})=\sigma(\widetilde{L}_{\beta})$ by Lemma \ref{lem:sigmainfini=sigma} \textbf{c)}.
    
    \item Let $L_0=\left(\mathrm{d}/\mathrm{d}z\right)^{\ell}$. Then we have $$\sigma(\widetilde{L}_{\beta}) \leqslant \ell + 2 \sigma(L_0) + \sigma(L_{\beta})$$ since $\sigma(z^{m-1} L_{\beta})=\sigma(L_{\beta})$. Here we use the fact (see \cite[Theorem 3, p. 16]{LepetitSize}) that  $$\forall L_1, L_2 \in \Qbar(z)\left[\gf{\mathrm{d}}{\mathrm{d}z}\right]\;, \quad \sigma(L_1 L_2) \leqslant \ord(L_1)+ 2 \sigma(L_1)+ \sigma(L_2).$$ 

Moreover, a basis of solutions of the equation $L_0(y(z))=0$ is $(1, z, \dots, z^{\ell-1})$, so that a fundamental matrix of solution of the system $y'=A_{L_0} y$ is the wronskian matrix $$Y=\begin{pmatrix} 1 & z & \dots & z^{\ell-1} \\ 0 & 1 & & (\ell-1) z^{\ell-2} \\ \vdots & 0 & \ddots & \vdots \\ 0 & 0 & \dots & (\ell-1)! \end{pmatrix},$$ that satisfies $Y^{(s)}=0$ for $s$ large enough. Hence $(A_{L_0})_s=Y_s Y^{-1}=0$ for $s$ large enough and $\sigma(L_0)=~0$.
    
    \item By applying Theorem \ref{th:chudnovskynilssongevreyandre} of Subsection \ref{subsec:chudnovskyeffectif} to the Nilsson-Gevrey series $z^{\beta} F(z)$, we obtain $$\sigma(L_{\beta}) \leqslant (1+\log(2)) \max\left(1, 2 \log(\db), \sigma(L) \right).$$ Finally, we have $$\sigma(\widetilde{L}_{\beta}) \leqslant \ell + (1+\log(2)) \max\left(1, 2 \log(\db), \sigma(L) \right).$$\end{itemize}
All this gives, by Theorem \ref{th:liensigmaYsigmaG}:
\begin{multline}\label{eq:defLambda0}
    \Lambda(\widetilde{L}_{\beta,\infty}) \leqslant (\mu+\ell)^2 \Big( \ell + (1+\log(2)) \max\left(1, 2 \log(\db), \sigma(L) \right)\Big)+(\mu+\ell)^2+\mu+\ell-1\\+(\mu+\ell)(\mu+\ell-1) H\left(\den(f_1,\dots,f_{\ell}) \db\right)+ ((\mu+\ell)^2+1) \sum\limits_{\zeta \neq 0, \infty} h\left(1/\zeta\right):=\Lambda_0(L,\beta),
\end{multline} where the last sum is on the roots of $\chi_L(z)=\sum\limits_{j=0}^{\ell} q_{j,\mu} z^j$, which doesn't depend on $\beta$. 
\end{dem}

\begin{rqu}
Since $\sigma(L) \leqslant \left(6 \mu^2(\delta+1)-1-(\mu - 1)(\delta+1)\right)\overline{\sigma}(F)$ by Theorem \ref{th:chudnovskyeffectif} above, we also have 
    $$\sigma(\widetilde{L}_{\beta}) \leqslant \ell + \big(1+\log(2)\big) \max\Big(1, 2 \log(\db),\big(6 \mu^2 (\delta+1)-1-(\mu - 1)(\delta+1)\big)\overline{\sigma}(F) \Big).$$
    This enables us to compute an upper bound on $\Lambda_0(L,\beta)$ in terms of $\overline{\sigma}(F)$ rather than $\sigma(L)$.
\end{rqu}

\subsection{Estimate of the denominator of $1/W_{\beta}(n)$}\label{subsec:estimatedenom1/Wbeta}

In order to estimate the denominator of $1/W_{\beta}(n)$, we are going to rely on a recurrence relation of order one satisfied by $(W_{\beta}(n))_{n \geqslant m}$ which enables us to express this sequence with Pochhammer symbols.

This is why we are going to need the following lemma. It provides estimates for the denominator of a quotient of Pochhammer symbols. It is also necessary for the proof of Lemma \ref{lem:denomCusn} in Subsection \ref{subsec:auxiliaryseries} above.

\begin{lem} \label{lem:denomquotientspochhammer}
Let $\alpha,\beta \in \Q \setminus \Z_{\leqslant 0}$. Then
\begin{enumerateth}
\item We have \begin{equation}\label{eq:quotientspochhammer1}\limsup_{n \rightarrow + \infty}  \left(\den\left(\gf{(\alpha)_0}{(\beta)_0},\dots,\gf{(\alpha)_n}{(\beta)_n}\right)\right)^{1/n} \leqslant \den(\alpha)^2 e^{\den(\beta)}.\end{equation}
\item In the case where $\beta=1$, for $n \in \N^*$, the binomial coefficient $\dbinom{\alpha+n-1}{n}:=\gf{(\alpha)_n}{n!}$ satisfies $$\den(\alpha)^{2n} \dbinom{\alpha+n-1}{n} \in \Z.$$ so that \begin{equation}\label{eq:quotientspochhammer2}\limsup_{n \rightarrow + \infty}  \left(\den\left(\gf{(\alpha)_0}{0!},\dots,\gf{(\alpha)_n}{n!}\right)\right)^{1/n} \leqslant \den(\alpha)^2.\end{equation}
\end{enumerateth}
\end{lem}

We give a proof for the reader's convenience. It is inspired by an argument of Siegel in \cite[pp. 56--57]{Siegel}. Note that \eqref{eq:quotientspochhammer2} saves a factor of $e$ from \eqref{eq:quotientspochhammer1}.

\begin{dem}
Write $\alpha=a/b$, $\beta=c/d$ and assume that $a$ and $b$ (resp. $c$ and $d$) are coprime. Let $$u_n=\gf{b^{2n}}{d^n} \gf{(\alpha)_n}{(\beta)_n}=\gf{b^n a(a+b)\dots(a+(n-1)b)}{c(c+d) \dots (c+(n-1)d)}.$$ If $p$ is a prime factor of $c(c+d)\dots(c+(n-1)d)$, then since $(c,d)=1$, $p$ does not divide $d$.  If $\ell \in \N^*$, then among $p^\ell$ consecutive integers of the form $c+\nu d$, $\nu \in \{0, \dots, n-1 \}$, only one is divisible by $p^\ell$. Thus at least $\left\lfloor n/p^\ell \right\rfloor$ and at most $\left\lfloor n/p^\ell \right\rfloor+1$ of the $c+\nu d$ are divisible par $p^\ell$. This proves that $$\sum\limits_{\ell=1}^{C_{n,p}} \left(\left\lfloor\gf{n}{p^\ell}\right\rfloor +1\right) \geqslant v_p(c(c+d) \dots (c+(n-1)d))\geqslant \sum\limits_{\ell=1}^{C_{n,p}} \left\lfloor\gf{n}{p^\ell}\right\rfloor=v_p(n!),$$ with $C_{n,p}=\left\lfloor \log(n)/\log(p)\right\rfloor$. Moreover, assuming that $p$ does not divide $b$, we get likewise $$\sum\limits_{\ell=1}^{C_{n,p}} \left(\left\lfloor\gf{n}{p^\ell}\right\rfloor +1\right) \geqslant v_p(a(a+b)\dots(a+(n-1)b))\geqslant \sum\limits_{\ell=1}^{C_{n,p}} \left\lfloor\gf{n}{p^\ell}\right\rfloor.$$ Hence $$v_p(u_n) \geqslant -C_{n,p}.$$

Notice that if $\beta=1$, we have $c=d=1$ and $v_p(c(c+d) \dots (c+(n-1)d))=v_p(n!)= \sum\limits_{\ell=1}^{C_{n,p}}\left\lfloor n/p^\ell \right\rfloor$, so that $v_p(u_n) \geqslant 0$, whence $\den(\alpha)^{2n} \gf{(\alpha)_n}{n!} \in \Z$. This proves \textbf{b)}.

\medskip

We come back to the general case $\beta \neq 1$. On the other hand, if $p$ divides $b$, then $$v_p(b^n a(a+b) \dots (a+(n-1)b)) \geqslant nv_p(b) \geqslant n,$$ and $$v_p(c(c+d) \dots (c+(n-1)d)) \leqslant v_p(n!)+C_{n,p}\leqslant n+C_{n,p}$$ so $v_p(u_n) \geqslant -C_{n,p}$. Therefore, $\Delta_n:=\prod\limits_{p \leqslant c+(n-1)d} p^{C_{n,p}}$  is a common denominator of $u_0, \dots, u_n$, so that $b^{2n} \Delta_n$ is a multiple of $\den\left( (\alpha)_0/(\beta)_0,\dots,(\alpha)_n/(\beta)_n\right)$.

And $$\log \Delta_n=\sum\limits_{p \leqslant c+(n-1)d} \log(p) \left\lfloor \gf{\log(n)}{\log(p)}\right\rfloor \leqslant \pi(c+(n-1)d) \log(n),$$ where $\pi$ is the prime-counting function. The Prime Number Theorem then implies that $$\pi(c+(n-1)d) \sim \gf{c+(n-1)d}{\log(c+(n-1)d)} \sim d \gf{n}{\log(n)},$$ so that $\log \Delta_n \leqslant d n+o(n)$. Since $b=\den(\alpha)$ and $d=\den(\beta)$, we thus obtain \textbf{a)}.\end{dem}

Let us now use this lemma to bound the denominator of $W_{\beta}(n)^{-1}$. We know that $Q_0(X)=\gamma_0 \prod\limits_{i=1}^{\mu} (X-e_i)$ and $Q_{\ell}(X)=\gamma_{\ell} \prod\limits_{i=1}^{\mu} (X+f_i-\ell)$, with $\gamma_0, \gamma_{\ell} \in \Oal_{\K}$. Thus, we find with the same computation as in \cite[pp. 13--14]{FRivoal} that $$\forall n \geqslant m \; , \quad Q_{\ell}(-n-\beta)W_{\beta}(n+1)=(-1)^{\ell} Q_0(-n-\beta)W_{\beta}(n)$$ and therefore that
\begin{equation}\label{eq:estimdenomWbeta-eq1}\forall n \geqslant m \; , \quad W_{\beta}(n)=W_{\beta}(m)\left((-1)^{\ell} \gf{\gamma_0}{\gamma_\ell}   \right)^{n-m} \prod_{i=1}^\mu  \frac{(m+ e_{i}+\beta )_{n-m}}{(m-f_{i} + \ell+\beta)_{n-m}}. \end{equation} Hence, by Lemma \ref{lem:denomquotientspochhammer}, \begin{equation}\label{eq:denomsection2.2}\limsup_{n \rightarrow + \infty} \den\left(\gf{1}{W_{\beta}(m)}, \gf{1}{W_{\beta}(m+1)}, \dots, \gf{1}{W_{\beta}(n)}\right)^{1/n} \leqslant \den(1/\gamma_0) \den(\mathbf{e}, \beta)^{2 \mu} \exp\left(\mu \den(\mathbf{f}, \beta)\right),\end{equation} where $\den(\mathbf{e},\beta):=\den(e_1, \dots, e_{\mu}, \beta)$ and $\den(\mathbf{f}, \beta):=\den(f_1,\dots,f_{\mu},\beta)$.

\subsection{Estimate of the denominator of the $D_{j,\beta}(n)/Q_{\ell}(1-n-\beta)$} \label{subsec:estimatedenomdjbetaql}

We consider $\mathcal{D}_n:=\den\left(u_{j,\beta}(k), 1 \leqslant j \leqslant \ell, m \leqslant k \leqslant n\right)$. Then by the determinant formula \eqref{eq:defDjbeta}, we have $$\forall m \leqslant k \leqslant n\; , \quad \mathcal{D}_{n+\ell-2}^{\ell-1} D_{j,\beta}(k) \in \Oal_{\K}.$$ Moreover,$$\gf{1}{Q_{\ell}(1-n-\beta)}=\gf{1}{\gamma_{\ell}} \prod_{i=1}^{\mu} \gf{1}{1-(n+\ell)+f_i-\beta}=\gf{\den(\mathbf{f},\beta)^{\mu}}{\gamma_{\ell}} \prod_{i=1}^{\mu} \gf{1}{\den(\mathbf{f}, \beta)(1-(n+\ell))+\nu_i}$$ with $\nu_i=\den(\mathbf{f},\beta)(f_i-\beta)$. For $n \geqslant m$, we have $$\forall k \leqslant n, \; \; \forall 1 \leqslant i \leqslant \mu, \quad \left\vert \den(\mathbf{f}, \beta)(1-(k+\ell)+\nu_i \right\vert \leqslant \den(\mathbf{f}, \beta)(n+\ell-1)+\nu_0,$$ where $\nu_0=\max |\nu_i|$. Thus, $$\forall n \geqslant m, \;\; \forall m \leqslant k \leqslant n, \quad d_{\den(\mathbf{f}, \beta)(n+\ell-1)+\nu_0}^{\mu} \den\left(1/\gamma_{\ell}\right) \gf{1}{Q_{\ell}(1-k-\beta)} \in \Oal_{\K},$$ with $d_s=\mathrm{lcm}(1,2, \dots, s)$. Finally, \begin{equation}\label{eq:denomsection2.3}
\forall m \leqslant k \leqslant n \;, \quad \mathcal{D}_{n+\ell-2}^{\ell-1} d_{\den(\mathbf{f}, \beta)(n+\ell-1)+\nu_0}^{\mu} \den\left(1/\gamma_{\ell}\right) \gf{D_{j,\beta}(n)}{Q_{\ell}(1-n-\beta)} \in \Oal_{\K}.\end{equation}

\bigskip

\subsection{Conclusion of the proof of Theorem \ref{th:calculC1etC2}} \label{subsec:computationC1Fbeta}

\subsubsection{Computation of $C_2(F,\beta)$}

Recall that $\delta_n$ is defined as the common denominator of the $1/W_{\beta}(k), D_{j,\beta}(k)/Q_{\ell}(1-k-\beta)$ and $u_{j,\beta}(k)$, when $k \in \{m, \dots, n \}$ and $j \in \{1, \dots, \ell \}$. Equations \eqref{eq:denomsection2.1}, \eqref{eq:denomsection2.2} and \eqref{eq:denomsection2.3} yield \begin{equation*}\begin{split}\limsup\limits_{n \rightarrow + \infty} \delta_n^{1/n} &\leqslant \exp\left(\max(1,\ell-1)[\K:\Q] \Lambda_0(L,\beta)\right) \exp\left(\den(\mathbf{f}, \beta)\right) \den\left(1/\gamma_0\right) \\ &\qquad\qquad\qquad\qquad\qquad\qquad\qquad\qquad\qquad \times \den(\mathbf{e}, \beta)^{2 \mu} \exp\left(\mu \den(\mathbf{f}, \beta)\right) 
\\
&\leqslant \den\left(1/\gamma_0\right) \den(\mathbf{e}, \beta)^{2 \mu} \exp\left(\max(1,\ell-1)[\K:\Q] \Lambda_0(L,\beta)+(\mu+1)\den(\mathbf{f}, \beta) \right),
\end{split}
\end{equation*}so that $$C_2(F,\beta):=\den\left(1/\gamma_0\right)^3 \den(\mathbf{e}, \beta)^{6 \mu} \exp\left(3\max(1,\ell-1)[\K:\Q] \Lambda_0(L,\beta)+3(\mu+1)\den(\mathbf{f}, \beta) \right)$$ is an upper bound on $\limsup \delta_n^{3/n}$.

As mentioned in the introduction, we observe that $C_2(F,\beta)$ only depends on $F$ and on the denominator $\db$ of $\beta$.

\bigskip

\subsubsection{Computation of $C_1(F)$}
In this part, we are going to find a suitable explicit constant $C_1(F,\beta)$ satisfying Proposition \ref{prop:recurrenceFns} \textbf{b)}. It will turn out not to depend on $\beta$.
\begin{itemize}
 \item Since any hypergeometric series of the form $\sum\limits_{k=0}^{\infty} \gf{(a^{(1})_k \dots (a^{(p)})_k}{(b^{(1})_k \dots (b^{(p)})_k} z^k$ has a radius of convergence~$1$, using formula \eqref{eq:estimdenomWbeta-eq1}, we have $$\limsup\limits_{n \rightarrow +\infty} \house{\gf{1}{W_{\beta}(n)}}^{1/n} \leqslant \house{\gf{\gamma_0}{\gamma_{\ell}}}.$$ 
    \item Let $j \in \{1, \dots, \ell \}$. We have seen in the proof of Proposition \ref{prop:denomujbeta} that if $u_n=u_{j,\beta}(n)$ and $U(z)=\sum\limits_{n=m}^{\infty} u_n z^n$, then $\widetilde{L}_{\beta,\infty}(U(z))=0$.

By Frobenius' Theorem (see \cite[Theorem 3.5.2, p. 349]{Hille}), the radius of convergence of $U$ around $0$ is exactly the largest $R >0$ such that the coefficients of $L_{\beta, \infty}$ are holomorphic in the punctured disc $D(0,R) \setminus \{ 0 \}$. Precisely, $R$ is the equal to the minimum of $|\xi|$, when $\xi$ is a nonapparent finite nonzero singularity of $\widetilde{L}_{\beta,\infty}$. Hadamard-Cauchy formula yields $$\limsup\limits_{n \rightarrow +\infty} |u_n|^{1/n} = \gf{1}{R}.$$ Furthermore, we have seen in Subsection \ref{subsec:estimatedenomujbeta} that the finite nonzero singularities of $\widetilde{L}_{\beta,\infty}$ are among the $\zeta^{-1}$, when $\zeta \neq 0$ is a root of $\chi_L(z)$. So $R \geqslant \min\limits_{\chi_L(\zeta)=0} |\zeta|^{-1}$, whence $$\limsup\limits_{n \rightarrow +\infty} |u_n|^{1/n} \leqslant \Phi_0(L),$$ with \begin{equation}\label{eq:definitionPhi0}
    \Phi_0(L):=\max_{\chi_L(\zeta)=0} |\zeta|\;,\;\; \text{where} \;\; \chi_L(z)=\sum_{j=0}^{\ell} q_{j,\mu}z^j \;\;\; \text{is independent of} \;\; \beta.
\end{equation} We recall that for all $j \in \{ 0, \dots, \ell \}$, $Q_j(X)=\sum\limits_{k=0}^{\mu} q_{j,k} X^j$.
\item The determinant formula \eqref{eq:defDjbeta} then implies that $$\limsup\limits_{n \rightarrow + \infty} \left|\gf{D_{j,\beta}(n)}{Q_{\ell}(1-n-\beta)}\right|^{1/n} \leqslant \Phi_0(L)^{\ell-1}.$$
\end{itemize}

Finally, $C_1(F):=\max\left(1,\house{\gamma_0/\gamma_{\ell}}, \Phi_0(L)^{\max(1,\ell-1)}\right)$ satisfies $$\max_{m+1 \leqslant k \leqslant n} \max\left(|u_{j,\beta}(k)|, \gf{1}{|W_{\beta}(k)|}, \gf{|D_{j,\beta}(k)|}{|Q_{\ell}(1-k-\beta)|} \right) \leqslant C_1(F)^{n(1+o(1))},$$ as claimed. This completes the proof of Theorem \ref{th:calculC1etC2}. Note that $C_1(F)$ does not depend on $\beta$.

\section{Examples}\label{sec:examples}

We now apply Theorem \ref{th:chapitre2} to some classical examples of $G$-functions, and compute explicitly $C(F,\beta)$ for them. The only real difficulty is to find an upper bound on $\sigma(L)$ where $L$ is the minimal operator of the $G$-function $F(z)$, but it can often be computed from the knowledge of the arithmetic behavior of the coefficients of $F(z)$. We are going to use the results of Section \ref{sec:sizeGop} to this end. We have implemented on Sage programs computing $C(F,\beta)$ assuming that we have a bound on $\sigma(L)$.

As mentioned in the preliminary remark of Section \ref{sec:explicitC}, the cases $\ell=1$ and $\ell \geqslant 2$ use different formulas -- the second case being much simpler. This is why the following examples are separated into these categories.

\subsection{Examples for which $\ell=1$}
\begin{itemize}
\item Let us first consider the simple example of $F(z)=\sum\limits_{k=0}^{\infty} z^k = 1/(1-z)$. The minimal operator of $F$ over $\Qbar(z)$ is $$L=(1-z)\left(\ddz\right)-1.$$ It satisfies $\ell=\delta=1$ and $\ell_0(0)=1$ (see \eqref{eq:defl0beta}). For any $n \in \N^*$ and $0 \leqslant s \leqslant S$, the $G$-function $F_{n,0}^{[s]}(z)$ is equal to \begin{equation}\label{eq:Fnbeta=polylogs} F_{n,0}^{[s]}(z)=\sum\limits_{k=0}^{\infty} \gf{z^{k+n}}{(k+n)^s}=\mathrm{Li}_s(z)-\sum\limits_{k=1}^{n-1} \gf{z^k}{k^s},\end{equation} so that the vector space spanned by the $F_{n,0}^{[s]}(\alpha)$, $n \in \N^*, 0 \leqslant s \leqslant S$ is equal to $\Vect\big(1, \mathrm{Li}_s(\alpha),\; 0 \leqslant s \leqslant S \big)$. For $\alpha \in \Qbar$, $0< |\alpha|<1$, Marcovecchio (\cite{Marcovecchio}) proved that $$\dim \Vect_{\Q(\alpha)}\big(1, \mathrm{Li}_s(\alpha),\; 0 \leqslant s \leqslant S \big) \geqslant \gf{1+o(1)}{[\Q(\alpha):\Q]\left(1+\log(2)\right)} \log(S),$$ thus obtaining a suitable constant $\widetilde{C}(F,\beta)=1+\log(2) \simeq 1.693$ for Theorem \eqref{th:chapitre2} in that case.

With our method, since $\ell=1$, we can use a refinement of \eqref{eq:expressionC2ell=1} allowed by Lemma \ref{lem:denomquotientspochhammer} \textbf{b)} and we find that $C(F,\beta)=4+\log(2) \simeq 4.693$. Thus, our method does not improve here the constant already known in the lower bound of Theorem \ref{th:chapitre2}.
\begin{rqu}
This operator $L$ is an example of "trivial case" for our estimates of Section \ref{sec:explicitC}, since $\ell=1$ and moreover $$\widetilde{L}_0=(1-z)\left(\ddz\right)^2-2\left(\ddz\right).$$ satisfies $\sigma(\widetilde{L}_0)=0$.
Indeed, let $A$ be the companion matrix of $\widetilde{L}_0$; then one can show by induction that $$A_s=\gf{(-1)^s s!}{(z-1)^s} \begin{pmatrix} 0 & 1-z \\ 0 & s+1 \end{pmatrix}.$$  
\end{rqu}

\medskip

If $\beta$ is an arbitrary positive rational number, the study of the sequence $\left(F_{n,\beta}^{[s]}(z)\right)_{n,s}$ provides a statement about the linear independence of the values of Lerch's functions, also quoted in \cite[p. 2]{Rivoal2003}. They are defined for $\mathrm{Re}(s) > 1$, $\beta \in \Q_{\geqslant 0}$ by $$\Phi_s(z,\beta)= \sum\limits_{k=0}^{\infty} \gf{z^k}{(k+\beta)^s}, \quad |z|<1.$$

Indeed, if $\alpha \in \Qbar$ and $0 < |\alpha|<1$, the same reasoning as in \eqref{eq:Fnbeta=polylogs} leads to $$\Vect_{\K}(F_{n,\beta}^{[s]}(\alpha), n \in \N, 0 \leqslant s \leqslant S)=\Vect_{\K}(1, \Phi_s(\alpha,\beta), 0 \leqslant s \leqslant S).$$ 

For example, if $\beta=1/2$, we numerically compute $C(F,\beta) \simeq 28.01$ so that at least $\gf{1+o(1)}{28.01} \log(S)$ of the values $\Phi_s(\alpha,1/2)$ for $0 \leqslant s \leqslant S$ are linearly independent over $\Q$ when $S \rightarrow +\infty$.
\bigskip

\item Let $\nu \geqslant 2$, $\boldsymbol{a}=(a_1, \dots, a_\nu) \in (\Q \setminus \Z_{\leqslant 0})^{\nu}$ and $\boldsymbol{b}=(b_1, \dots, b_{\nu-1}) \in (\Q \setminus \Z_{\leqslant 0})^{\nu-1}$. The hypergeometric function $F(z)=_\nu F_{\nu-1}(\boldsymbol{a}, \boldsymbol{b} ; z)$ defined in the introduction as $$_{\nu} F_{\nu-1}(\boldsymbol{a}, \boldsymbol{b} ; z)=\sum\limits_{k=0}^{\infty} \gf{(a_1)_k \dots (a_{\nu})_k}{k! (b_1)_k \dots (b_{\nu-1})_k} z^k$$ is a nonpolynomial $G$-function that is solution of the generalized hypergeometric equation $\mathcal{H}_{\boldsymbol{a}, \boldsymbol{b}}(y(z))=0$ where $$\mathcal{H}_{\boldsymbol{a}, \boldsymbol{b}}=z(\theta+a_1) \dots (\theta+a_{\nu})-(\theta+b_1-1) \dots (\theta+b_{\nu-1}-1)\theta \in \Q\left[z,\ddz\right].$$

If for all $i,j$, $a_i - b_j \not\in \Z$, then $\mathcal{H}_{\boldsymbol{a},\boldsymbol{b}}$ is a minimal operator over $\Qbar(z)$ for $F(z)$. It satisfies $\mu=\nu$, $\delta \leqslant \nu+1$ and $\omega=\nu$, so $\ell=1$. 

Hence, by Theorem \ref{th:chudnovskyeffectif}, \begin{align*}\sigma(\mathcal{H}_{\boldsymbol{a},\boldsymbol{b}}) &\leqslant \left(5 \nu^2(\nu+2)-1-(\nu-1)(\nu+2)\right)\overline{\sigma}(F) \\ 
&\leqslant \left(5 \nu^2(\nu+2)-1-(\nu-1)(\nu+2)\right)\big(2 \nu \log(\den(\boldsymbol{a})+(\nu-1)\den(\boldsymbol{b})\big) \end{align*} since $F(z)$ has radius of convergence $1$ and using the estimate on the denominator of a quotient of Pochhammer symbols given by Lemma \ref{lem:denomquotientspochhammer}.

The exponents of $\mathcal{H}_{\boldsymbol{a},\boldsymbol{b}}$ at $0$ are $0$, $1-b_1$, ..., $1-b_{\nu-1}$ and its exponents at $\infty$ are $a_1$, ..., $a_{\nu}$.

For example, for $\nu=2$, $a_1=1/3$, $a_2=2/11$, $b_1=1/6$, and $\beta=\gf{1}{7}$, we find numerically that, on the one hand $\sigma\left(\mathcal{H}_{\mathbf{a}, \mathbf{b}}\right) \leqslant 1119 $ and on the other hand $C(F,\beta) \simeq 2443$ so that, if $\alpha \in \Qbar^*$ and $0<|\alpha|<1$, at least  $\gf{1+o(1)}{2443 [\Q(\alpha):\Q]} \log(S)$ of the numbers $$\sum\limits_{k=0}^{\infty} \gf{(1/3)_k (2/11)_k}{k! (1/6)_k} \gf{z^k}{(7k+1)^s}, \quad 0 \leqslant s \leqslant S$$ are linearly independent when $S$ tends to infinity. This provides a refinement of \cite[Corollary 1, p. 4]{FRivoal} in that case, where the constant was not given.
\end{itemize}

\subsection{Examples for which $\ell \geqslant 2$}
\begin{itemize}
\item We first consider the following $G$-function which is algebraic of $\Qbar(z)$: $$F(z)=\dfrac{1}{\sqrt{1-6z+z^2}}=\sum_{k=0}^{\infty} \sum\limits_{j=0}^{k} \binom{k+j}{j} \binom{k}{j} z^k$$ whose minimal operator over $\Qbar(z)$ is
$$L=(z^2-6z+1)\left(\ddz \right)+z-3.$$
The power series $F(z) \in \Z\llbracket z \rrbracket$ has radius of convergence $3-2\sqrt{2}$, hence $\overline{\sigma}(F) \leqslant -\log(3-2\sqrt{2})$. Furthermore, $L$ satisfies $\ell=\delta=2$, $\mu=1$. Its exponent at $0$ (resp. $\infty$) is $0$ (resp. $1$).

With $\beta=3/5$, we obtain $\ell_0(\beta)=2$ and $C(F,\beta) \simeq 3164$, so that, for every $\alpha \in \Qbar^*$, $0<|\alpha|<3-2\sqrt{2}$, at least $\gf{1+o(1)}{3164[\Q(\alpha):\Q]} \log(S)$ of the numbers $$\sum_{k=0}^\infty\sum_{j=0}^k \binom{k}{j}\binom{k+j}{j} \frac{\alpha^k}{(5k+3)^s} \quad \textup{and} \quad \sum_{k=0}^\infty \sum_{j=0}^k \binom{k}{j}\binom{k+j}{j} \frac{\alpha^k}{(5k+8)^s}, \quad 0 \leqslant s \leqslant S.
$$
are linearly independent over $\Q(\alpha)$ when $S \rightarrow +\infty$.


\medskip
For $\beta=0$, we compute $\ell_0(0)=2$ and $C(F,0) \simeq 3083$.

\bigskip

\item Denote, for all $k \in \N^*$, $\mathcal{L}_k=\left((1-z)\ddz-1\right)^k \times \ddz$ the minimal operator of $\log(1-z)^k$. It is an operator of order $k+1$, with $\omega=0$ and $\delta=k$. Its exponents at $0$ (resp. $\infty$) are $0, 1, \dots, k$ (resp. $0$).

Using \cite[Proposition 6, p. 12]{LepetitSize}, we obtain that $\sigma(\mathcal{L}_k) \leqslant (1+\log(k))$.

For $k=2$ and $\beta=1/12$, we have $\ell_0(\beta)=2$ and we obtain a constant $C(F,\beta) \simeq 1912$.

Since $\log(1-z)^2=2\sum\limits_{k=0}^{\infty} \gf{1}{k+1} \left(\sum\limits_{i=1}^{k} \gf{1}{i}\right) z^{k+1}$, this yields for $\alpha \in \Qbar^*$, $0<|\alpha|<1$, 
\begin{multline*}\dim_{\Q(\alpha)}\mathrm{Span}\left( \sum\limits_{k=0}^{\infty} \gf{1}{k+1} \left(\sum\limits_{i=1}^{k} \gf{1}{i}\right) \gf{\alpha^{k+1}}{(12k+1)^s}\;,\quad \sum\limits_{k=0}^{\infty} \gf{1}{k+1} \left(\sum\limits_{i=1}^{k} \gf{1}{i}\right) \gf{\alpha^{k+1}}{(12k+13)^s}, \quad 0 \leqslant s \leqslant S \right) \\ \geqslant \gf{1+o(1)}{1912 [\Q(\alpha):\Q]} \log(S).\end{multline*}


For $\beta=0$, we compute $\ell_0(0)=2$ and $C(F,0) \simeq 570$.

\bigskip

\item We finally consider the example of the generating series of the Apéry numbers, introduced by Apéry in \cite{Apery} in order to prove the irrationality of $\zeta(3)$: $$F(z)=\sum_{k=0}^{\infty} \sum_{j=0}^k \dbinom{k+j}{j}^2 \dbinom{k}{j}^2 z^k \in \Z\llbracket z \rrbracket$$
It is a $G$-function of radius of convergence $(\sqrt{2}-1)^4 \leqslant 1$ and its coefficients are integers so $\overline{\sigma}(F) \leqslant -4 \log(\sqrt{2}-1)$. The minimal nonzero differential operator of $F(z)$ over $\Qbar(z)$ is $$L=z^2(1-34z+z^2) \left(\ddz\right)^3+z(3-153z+6z^2)\left(\ddz\right)^2+(1-112z+7z^2)\ddz+z-5.$$
It is of order $3$, with $\delta=4$ and $\ell=2$. Using Theorem \ref{th:chudnovskyeffectif}, we obtain, with $\beta=2/3$, $C(F,\beta) \simeq 209532.$

Moreover, $\ell_0(\beta)=2$. Thus, if $\alpha \in \Qbar^*$ and $0 < |\alpha| < \left(\sqrt{2}-1\right)^4$, the dimension over $\Q(\alpha)$ of the vector space generated by the numbers $$ \sum_{k=0}^\infty \sum_{j=0}^k \binom{k}{j}^2\binom{k+j}{j}^2 \frac{\alpha^k}{(3k+5)^s} \quad \textup{and} \quad \sum_{k=0}^\infty \sum_{j=0}^k \binom{k}{j}^2\binom{k+j}{j}^2 \frac{\alpha^k}{(3k+8)^s}, \quad 0 \leqslant s \leqslant S.$$
is at least $\gf{1+o(1)}{209532 [\Q(\alpha):\Q]} \log(S)$ when $S \rightarrow +\infty$.

\medskip
For $\beta=0$, we find $\ell_0(0)=2$ and $C(F,0) \simeq 209533$.
\end{itemize}

\printbibliography

\bigskip
G. Lepetit, Université Grenoble Alpes, CNRS, Institut Fourier, 38000 Grenoble, France.

\url{gabriel.lepetit@univ-grenoble-alpes.fr}. 

\bigskip

\emph{Keywords} : $G$-functions, $G$-operators, Linear independence criterion, Saddle point method.

\bigskip

\emph{2020 Mathematics Subject Classification}. Primary 11J72, 11J91 Secondary 34M03, 34M35, 41A60.

\end{document}